\documentclass[a4paper,11pt]{article}
\pdfoutput=1

\usepackage[subsec]{encstyle}

\usepackage{fullpage}

\title{Higher category theory}
\author{Rune Haugseng}

\date{\today}

\renewcommand{\Top}{\catname{Top}}

\begin{document}

\maketitle
\begin{abstract}
  Invited contribution to the \emph{Encyclopedia of Mathematical
    Physics}. We give an introduction to the homotopical theory of
  higher categories, focused on motivating the definitions of the
  basic objects, namely $\infty$-categories and
  $(\infty,n)$-categories.
\end{abstract}

\tableofcontents

\section{What are higher categories?}
Higher categories have recently become an important tool in several
areas of mathematics. In this article we give a brief
introduction to the modern, homotopical theory of higher categories,
focusing on describing the basic objects of interest and explaining
how they are defined. In this section we start by giving a first idea
of what a higher category should be in \S\ref{subsec:ncats} and
informally describing some examples in \S\ref{subsec:exs}. We then
explain why it is non-trivial to give a useful definition of higher
categories in \S\ref{subsec:weak}.  In \S\ref{subsec:hohyp} we
introduce the homotopical approach to higher categories by discussing
Grothendieck's Homotopy Hypothesis, after giving some topological
background in \S\ref{subsec:htpy}.

\subsection{An informal description of $n$-categories}\label{subsec:ncats}

The basic idea of a higher category is that it should be a structure
that has, in addition to the objects and morphisms of an ordinary
category, additional layers of ``higher'' morphisms. Thus an
\emph{$n$-category} is a structure where we have:
\begin{itemize}
\item objects\footnote{It is sometimes convenient to think of objects
    as \emph{$0$-morphisms}.} ($\bullet$),
\item morphisms (or $1$-morphisms) between objects ($\bullet \to \bullet$),
\item $2$-morphisms between morphisms (with the same source and target),
  which we can depict as:
  
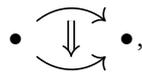
\[
    \begin{tikzcd}
      \bullet \ar[r, bend left=40, ""{below,name=A,inner sep=2pt}]
      \ar[r, bend right=40, ""{above,name=B,inner sep=2pt}] & \bullet
      \ar[from=A, to=B, Rightarrow]
    \end{tikzcd}
   \!\!,    
  \]
\item $3$-morphisms between $2$-morphisms, which we can depict as:
    
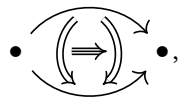
\[
    \begin{tikzcd}[column sep=large]
      \bullet \ar[r, bend left=50, ""{below,name=A,inner sep=2pt}]
      \ar[r, shorten <=13pt, shorten >=13pt, triple]
      \ar[r, bend right=50, ""{above,name=B,inner sep=2pt}] & \bullet
      \ar[from=A, to=B, bend right=50, shift right, Rightarrow, ""{below,name=C,
        inner sep=2pt}]      
      \ar[from=A, to=B, bend left=50, shift left, Rightarrow, ""{above,name=D,
        inner sep=2pt}]
    \end{tikzcd}
   \!\!,    
  \]
\item $4$-morphisms between $3$-morphisms, which we can depict as:
    
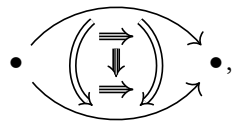
\[
    \begin{tikzcd}[column sep=small]
      \bullet \ar[rrr, bend left=50, ""{below,name=A}]
      \ar[rrr, bend right=50, ""{above,name=B}]
      & {} \ar[r, triple, shift left=2ex, ""{below,name=C,inner
        sep=2pt}] \ar[r, triple, shift right=2ex, ""{above,name=D,inner
        sep=2pt}]  & {} & \bullet
     \ar[from=A, to=B, bend right=50, shift right=1.2ex, Rightarrow]
     \ar[from=A, to=B, bend left=50, shift left=1.2ex, Rightarrow]
     \ar[from=C, to=D, quadruple]
   \end{tikzcd}
   \!\!,
  \]
\item   \ldots,
\item $n$-morphisms between $(n-1)$-morphisms.
\end{itemize}
We should be able to compose $i$-morphisms for all $0 < i \leq n$, and
an $i$-morphism should have an identity $(i+1)$-morphism for
$0 \leq i < n$. In particular, a $0$-category is just a set and a
$1$-category is an ordinary category.

An $i$-morphism in an $n$-category has a unique source and target
$j$-morphism for $j = 0,1,\dots,i-1$, and we should in fact be able to
compose two $i$-morphisms whose source and target $j$-morphisms agree
for each choice of $j$ --- for example, we can combine $2$-morphisms both
vertically and horizontally:
  
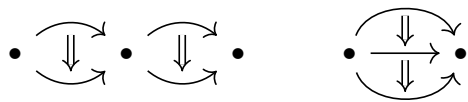
\[
    \begin{tikzcd}
      \bullet \ar[r, bend left=40, ""{below,name=A,inner sep=2pt}]
      \ar[r, bend right=40, ""{above,name=B,inner sep=2pt}] & \bullet
      \ar[r, bend left=40, ""{below,name=C,inner sep=2pt}]
      \ar[r, bend right=40, ""{above,name=D,inner sep=2pt}] & \bullet
      \ar[from=A, to=B, Rightarrow]
      \ar[from=C, to=D, Rightarrow]
      &  \bullet \ar[r, bend left=70, ""{below,name=E,inner sep=1pt}]
      \arrow[r, ""{above,name=F,inner sep=1pt}, ""{below,name=G,inner
        sep=1pt}]
      \ar[r, bend right=70, ""{above,name=H,inner sep=1pt}]
      & \bullet
      \ar[from=E, to=F, Rightarrow]
      \ar[from=G, to=H, Rightarrow]      
    \end{tikzcd}
  \]
These compositions must be compatible, so that there is for instance a
unique way to compose a diagram of $1$- and $2$-morphisms such as
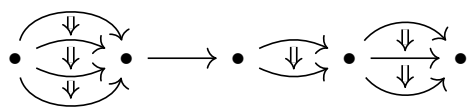
\[
\begin{tikzcd}
  \bullet 
  \ar[r, bend left=75, ""{below,name=A,inner sep=0.5pt}] 
  \ar[r, bend left=25, ""{name=B1,inner sep=0.5pt},
  ""{name=B2,below,inner sep=1pt}] 
  \ar[r, bend right=25, ""{name=C1,inner sep=1pt}, ""{below,name=C2,inner sep=0.5pt}] 
  \ar[r, bend right=75, ""{name=D,inner sep=1pt}] &
  \bullet 
  \ar[r]  &
  \bullet 
  \ar[r, bend left=25, ""{below,name=E,inner sep=1pt}] 
  \ar[r, bend right=25, ""{name=F,inner sep=1pt}] &
  \bullet 
  \ar[r, bend left=50, ""{below,name=G,inner sep=1pt}] 
  \ar[r, ""{name=H1,inner sep=1pt}, ""{name=H2,below,inner sep=1pt}] 
  \ar[r, bend right=50, ""{name=I,inner sep=1pt}] &
  \bullet
  \arrow[from=A,to=B1,Rightarrow]
  \arrow[from=B2,to=C1,Rightarrow]
  \arrow[from=C2,to=D,Rightarrow]
  \arrow[from=E,to=F,Rightarrow]
  \arrow[from=G,to=H1,Rightarrow]
  \arrow[from=H2,to=I,Rightarrow]
\end{tikzcd}
\]  
to a single $2$-morphism.

Slightly more precisely, if $\eC$ is an $n$-category, then for any pair of objects $x,y$ we
have an $(n-1)$-category $\eC(x,y)$ whose objects are morphisms from
$x$ to $y$, with morphisms being $2$-morphisms among these, etc., and we
can think of composition (with compatible source and target objects) as giving us functors of $(n-1)$-categories
\[ \eC(x,y) \times \eC(y,z)\to \eC(x,z);\] vertical composition of $2$-morphisms is then the composition of
$1$-morphisms in these $(n-1)$-categories, and so forth.

\subsection{Some examples}\label{subsec:exs}

At first glance, the idea of an $n$-category might seem just as
innocuous for general $n$ as it is in the case $n=1$, but it is in
fact not so easy to give a good definition. Before we explain why, let
us informally describe some structures that ought to give examples of
$n$-categories:\footnote{We revisit these examples in \cref{rmk:exs} at the end of our
  discussion of $(\infty,n)$-categories, and we postpone references to
actual definitions until then.}

\begin{ex}\label{ex:catofncat}
  The prototypical example of a category is the category of sets.  In
  the same way, the prototypical example of a $2$-category is the
  $2$-category of categories. This has (small) categories as objects,
  functors as $1$-morphisms, and natural transformations as
  $2$-morphisms. More generally, we expect to have an $(n+1)$-category
  of $n$-categories: Recall that a natural transformation between
  functors from $\uC$ to $\uD$ can be
  defined as a functor $\uC \times C_{1} \to \uD$, where $C_{1}$ is
  the ``universal morphism''.
  Similarly, an $i$-morphism between $n$-categories should be a
  functor $\eC \times C_{i} \to \eD$ where $C_{i}$ is the ``universal
  $i$-morphism''; informally, $C_{i}$ has two objects $0,1$ with
  \[
    C_{i}(s,t) =
    \begin{cases}
      \{\id_{s}\}, & s = t, \\
      \emptyset, & s = 1, t = 0, \\
      C_{i-1}, & s = 0, t = 1,
    \end{cases}
  \]
  so that we have
  \[ C_{0}=\bullet,\quad C_{1} = \bullet \to \bullet, \quad C_{2} =
        \begin{tikzcd}
      \bullet \ar[r, bend left=40, ""{below,name=A,inner sep=2pt}]
      \ar[r, bend right=40, ""{above,name=B,inner sep=2pt}] & \bullet
      \ar[from=A, to=B, Rightarrow]
    \end{tikzcd}\!\!,
  \]
  etc.
\end{ex}

\begin{ex}\label{ex:morita}
  Let $R$ be a commutative ring. Then we can define the \emph{Morita
    $2$-category} $\catname{Mor}(R)$ of $R$. This has associative $R$-algebras as objects,
  a $1$-morphism from $A$ to $B$ is an $A$-$B$-bimodule, and the
  $2$-morphisms are bimodule homomorphisms. The composite of an
  $A$-$B$-bimodule $M$ and a $B$-$C$-bimodule $N$ is the relative
  tensor product $M \otimes_{B} N$.
\end{ex}

\begin{ex}\label{ex:span}
  Let $\uC$ be a category with pullbacks. Then we should have a
  \emph{span $n$-category} $\Span_{n}(\uC)$ of $\uC$: This has the
  objects of $\uC$ as its objects, but its morphisms from $x$ to $y$
  are \emph{spans} (or \emph{correspondences}), that is diagrams
  \[ x \from z \to y \]
  in $\uC$; we compose two spans by taking pullbacks, so that in the
  following diagram the
  composite of the two red spans is the outer blue span:
  
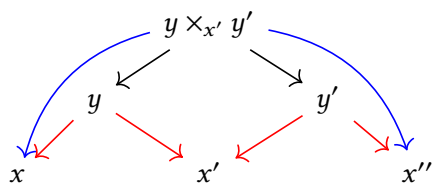
\[
    \begin{tikzcd}[column sep=small, row sep=small]
        &   & y \times_{x'} y' \ar[dl] \ar[dr] \ar[ddll, color=blue,
        bend right] \ar[ddrr, color=blue, bend left] \\
        & y \ar[dl,color=red] \ar[dr,color=red] & & y'
        \ar[dl,color=red] \ar[dr,color=red] \\
      x & & x' & & x''
    \end{tikzcd}
  \]
  Next, the $2$-morphisms in $\Span_{n}(\uC)$ are ``spans of spans'',
  that is diagrams of the form
  
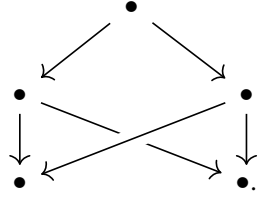
\[
    \begin{tikzcd}
       & \bullet \ar[dl] \ar[dr] \\
      \bullet \ar[d] \ar[drr] & & \bullet \ar[d]  \ar[dll, crossing over] \\
      \bullet & & \bullet.
    \end{tikzcd}
  \]
  Composition is again by taking pullbacks, and we keep considering
  iterated spans to define the $i$-morphisms; we can also say that the
  mapping $(n-1)$-category $\Span_{n}(\uC)(x,y)$ is
  $\Span_{n}(\uC_{/x,y})$, where $\uC_{/x,y}$ denotes the category of
  objects of $\uC$ equipped with morphisms to both $x$ and $y$.
\end{ex}

\begin{ex}\label{ex:cob}
  If $M$ and $N$ are two closed $k$-manifolds, a \emph{cobordism}
  between $M$ and $N$ is a compact $(k+1)$-manifold with boundary $X$
  together with a diffeomorphism $\partial X \cong M \amalg N$. We can
  define the \emph{cobordism category} $\catname{Cob}_{k,k+1}$ whose
  objects are closed $(k-1)$-manifolds and whose morphisms are
  $k$-dimensional cobordisms; composition is given by gluing
  cobordisms along their common boundary component. By considering
  manifolds with corners of higher codimensions we can extend this to
  define an $n$-category $\catname{Cob}_{k,k+n}$ where the $2$-morphisms
  are $(k+2)$-dimensional manifolds whose boundary decomposes
  appropriately into two $(k+1)$-dimensional cobordisms between the
  same closed $k$-manifolds. In particular, there is an $n$-category
  $\catname{Cob}_{0,n}$ whose $i$-morphisms are
  $i$-dimensional cobordisms with corners. (One can also consider
  variants where all the manifolds involved are equipped with
  compatible additional structures, such as orientations.)
\end{ex}

\begin{remark}
  Although we do not have space to discuss applications of higher
  categories in physics in this article, we should at least mention
  that most of the examples we have just outlined are in fact highly
  relevant to mathematical physics: The cobordism $n$-category
  $\catname{Cob}_{0,n}$ is central to the mathematical formalization
  of (extended) topological quantum field theories (TQFTs), and the
  Morita $2$-category and its higher-dimensional cousins are relevant as
  targets for certain interesting TQFTs (in particular those defined
  by \emph{factorization homology}); on the other hand, higher
  categories of spans are relevant when describing the structure of
  \emph{classical} field theories. In fact, TQFTs were one of the main
  early motivations for the development of higher category theory. We
  refer to other articles in this volume for detailed discussion of
  these topics.
\end{remark}

\subsection{Strict and weak $n$-categories}\label{subsec:weak}
We now want to explain why it is not so easy to give a ``correct''
definition of $n$-categories. To start, we note that an important
feature of $n$-categories is that we get, by induction, an
increasingly refined
notion of when two objects in an $n$-category are ``the same'', and
more generally of when an $i$-morphism is ``invertible'':
\begin{itemize}
\item In a set, two objects $x,y$ are the same if they are
  \emph{equal}: $x = y$
\item In a category $\eC$, two objects $x, y$ are the same if they are
  \emph{isomorphic}: there are morphisms $f \colon x \to y$, $g \colon
  y \to x$ such that $gf = \id_{y}$ in the set $\eC(y,y)$ and $fg =
  \id_{x}$ in the set $\eC(x,x)$; in this case we also say that $f$ (and $g$)
  are \emph{invertible} morphisms or \emph{isomorphisms} in $\eC$.
\item In a $2$-category $\eC$, two objects $x, y$ are the same if they are
  \emph{equivalent}: there are morphisms $f \colon x \to y$, $g \colon
  y \to x$ such that $gf$ is \emph{isomorphic} to $\id_{y}$ in the
  category $\eC(y,y)$ and $fg$ is isomorphic to $\id_{x}$ in the
  category $\eC(x,x)$; here we also say that $f$ (and $g$)
  are \emph{invertible} morphisms or \emph{equivalences} in $\eC$. We
  also define the invertible $2$-morphisms to be those that are
  isomorphisms in the mapping categories $\eC(x,y)$.
\item In an $n$-category $\eC$, two objects $x, y$ are the same if they are
  \emph{equivalent}: there are morphisms $f \colon x \to y$, $g \colon
  y \to x$ such that $gf$ is \emph{equivalent} to $\id_{y}$ in the
  $(n-1)$-category $\eC(y,y)$ and $fg$ is equivalent to $\id_{x}$ in the
  $(n-1)$-category $\eC(x,x)$; here we again say that $f$ (and $g$)
  are \emph{invertible} morphisms or \emph{equivalences} in $\eC$. We
  also say that the invertible $i$-morphisms for $i > 1$ are those
  that give invertible $(i-1)$-morphisms in the mapping
  $(n-1)$-categories $\eC(x,y)$.
\end{itemize}

In category theory, a key insight is that we should never ask for two
objects of a category to be equal, only that they are isomorphic under
a specified isomorphism. For example, the right notion of two functors
being ``the same'' is that they are \emph{naturally isomorphic}, not
equal, and therefore the right notion of two categories $\uC$ and
$\uD$ being ``the same'' is not that they are isomorphic, but that
there exist functors $F \colon \uC \to \uD$ and $G \colon \uD \to \uC$
with natural isomorphisms $GF \cong \id_{\uC}$, $FG \cong
\id_{uD}$.\footnote{Note that this is precisely the definition of $\uC$ and
$\uD$ being equivalent in the $2$-category of categories, functors, and
natural transformations.} It is reasonable to expect the analogous
principle to apply for $n$-categories: we should never ask for two
objects of an $n$-category to be equal, or even isomorphic, but only
\emph{equivalent} in the sense we just discussed.

If we accept this principle, we can appreciate the key difficulty in
defining $n$-categories ``correctly'': At first it may seem perfectly
reasonable to demand that the composition of $i$-morphisms should be
strictly associative, \ie{} there should be associativity identities
\[f(gh) = (fg)h\] for all composable $i$-morphisms; this gives the
notion of \emph{strict} $n$-categories.\footnote{These are also easy to
define inductively, for example a strict $n$-category is precisely a
category \emph{enriched} in strict $(n-1)$-categories; see \S\ref{subsec:enr}.}
However, here we are asking for two objects in an $(n-i)$-category to be \emph{equal}
rather than equivalent --- our principle tells us that we should instead 
supply an invertible
$(i+1)$-morphism $f(gh) \to (fg)h$ (the ``associator''). Now using these,
there are two ways to relate the different orders in which we may
compose 4 $i$-morphisms:
\begin{equation}
  \label{eq:pentagon}
  \begin{tikzcd}[column sep=tiny]
    {} &  & (fg)(hk) \arrow{drr}\\
    f(g(hk)) \arrow{dr}\arrow{urr} & & & & ((fg)h)k \\
    & f((gh)k) \arrow{rr} & & (f(gh))k. \arrow{ur}
  \end{tikzcd}
\end{equation}
To get a good notion of assocativity, these two composites
$f(g(hk)) \to ((fg)h)k$ should be ``the same'', which means they ought
to be related by a (specified) invertible $(i+2)$-morphism; using
these we can in turn formulate a coherence condition for composites of
5 morphisms in terms of an invertible $(i+3)$-morphism, and so on for
ever\footnote{The shapes of the coherence diagrams for compositions of
  increasing length are precisely the so-called \emph{associahedra},
  which are a family of polyhedra first introduced by
  Stasheff~\cite{Stasheff} to describe multiplications on topological
  spaces that are associative up to a coherent choice of higher
  homotopies.}
(or at least until we reach the
$n$-morphisms). This is the idea of \emph{weak} $n$-categories. Since
this coherence data quickly becomes intractable to write out
explicitly, we might think that the strict definition is the better
one. Unfortunately, almost all interesting examples of $n$-categories
fail to be strict --- this is true even in the simplest possible case:
\begin{ex}
  A $2$-category $\eC$ with a single object $*$ amounts to the data of a
  category $\uC = \eC(*, *)$, equipped with a functor
  $\blank\otimes\blank \colon \uC \times \uC \to \uC$ (composition of
  endomorphisms of $*$) and an object $\bbone \in \uC$ (the identity
  morphism of $*$). We can ask for the ``multiplication'' $\otimes$ to
  be strictly associative and unital, but such structures are
  exceedingly rare --- in practice, the useful notion of a tensor
  product on a category is that of a \emph{monoidal structure}, where
  we instead ask for \emph{natural isomorphisms}
  \[ X \otimes (Y \otimes Z) \cong (X \otimes Y) \otimes Z,\]
  \[ \bbone \otimes X \cong X \cong X \otimes \bbone,\] which must
  further satisfy certain coherence conditions; in particular, for the associativity
  isomorphisms 
  the two compositions in the pentagon
  \cref{eq:pentagon} must be equal for any quadruple tensor product.
\end{ex}
More generally, almost all of the examples we mentioned earlier can
only be defined as \emph{weak} $n$-categories.\footnote{The only
  exception is that there is a strict $(n+1)$-category of strict
  $n$-categories.} For $n=2$, it turns out that any weak $2$-category
is equivalent to a strict one\footnote{The earliest explicit reference
  I have found for this result is \cite{MacLanePare}, but the authors
  there say it is a special case of \cite{Ben68}*{Th\'eor\`eme
    5.2.4}}, so when working with $2$-categories we can in a sense
``get away with'' the strict theory, but this is false\footnote{See
  \cite{GordonPowerStreet}*{\S 8.5} --- there it is shown that weak
  $3$-categories with one object and one $1$-morphism are precisely
  braided monoidal categories \cite{JoyalStreet}. On the other hand,
  the Eckmann--Hilton argument shows that a strict $3$-category with
  one object and one $1$-morphism is a strict symmetric monoidal
  category, which is a far more restrictive structure.} for
$n \geq 3$. To get a theory of $n$-categories that encompasses the
examples we are interested in, we therefore have no choice but to
consider the weak version.

For $n = 2$, it is not too hard to write out an explicit definition of
a weak $2$-category, and there is an extensive literature on both strict
and weak $2$-categories (often called \emph{bicategories}); both were
originally defined by B\'enabou, in \cite{Ben65} and \cite{Ben67},
respectively. There is also some work on weak $3$-categories (or
\emph{tricategories}), with the first definition due to Gordon, Power
and Street~\cite{GordonPowerStreet}, but the coherence data for
composition in an $n$-category quickly becomes impossible to write out
explicitly\footnote{The unconvinced reader may peruse Trimble's
  explicit definition of a weak $4$-category \cite{Trimble}}. However,
it is possible to give systematic descriptions of the coherence data
and thus obtain definitions of weak $n$-categories for general $n$.
Various definitions of this type\footnote{The first definition of weak
  $n$-categories for arbitrary $n$, which was proposed by Street in
  \cite{StreetOriented}, has a somewhat different flavour; see
  \cref{rmk:otherinftyn}.}  were the focus of work on higher
categories in the 1990s and early 2000s; this includes in particular the definitions
of Baez and Dolan~\cite{BDopetope}, Batanin~\cite{Batanin}, and Leinster~\cite{LeinsterOpds}, which were perhaps the most
prominent.\footnote{All three of these definitions make use of various
  types of generalized \emph{operads} to describe the coherence
  data. Operads are objects that describe different
  types of algebraic structures; they were first introduced by
  May~\cite{May} and Boardman--Vogt~\cite{BoardmanVogt} to describe
  multiplications on topological spaces that are homotopy-coherently
  commutative.}

Such definitions have turned out to be difficult to work with in
practice, however, and often do not lend themselves easily to defining
interesting examples of higher categories. A key insight for work on
higher categories over the last two decades has been that instead of a
``bottom-up'' approach, where we try to combinatorially define
$n$-categories for increasing $n$, it is much easier to start by
defining \emph{$(\infty,1)$-categories} and then use them to define
more general higher categories. To explain what these objects are, we
first need the following terminology:

\begin{defn}
  Among the $n$-categories, we can single out those where all
  $i$-morphisms are invertible for $i > k$; these are called
  \emph{$(n,k)$-categories}. In the extreme cases we recover all
  $n$-categories as $(n,n)$-categories, while for $k = 0$ we get those
  $n$-categories where all $i$-morphisms are invertible for all
  $i \leq n$; these are also called
  \emph{$n$-groupoids}\footnote{Recall that a \emph{groupoid} is a
    category where all morphisms are isomorphisms.}.
\end{defn}

Informally, we can imagine a version of higher categories where we
keep going for ever and add $i$-morphisms for all $i$ instead of
stopping at some fixed $n$; such objects are called
\emph{$\omega$-categories} or \emph{$(\infty,\infty)$-categories}. If
we assume that the $i$-morphisms are in fact all
invertible\footnote{Although our discussion here is very informal, to
  avoid confusion let us mention that there are problems with
  defining what it means for a morphism in an $\omega$-category to be
  ``invertible'' in the bottom-up perspective we have so far
  considered; see \cref{rmk:coind}.}  for $i > k$, we get the notion
of \emph{$(\infty,k)$-categories}. For $k = 0$ these are also called
\emph{$\infty$-groupoids} and for $k = 1$ we will call them
\emph{$\infty$-categories} (following \cite{HTT} and most of the
subsequent literature).

\begin{warning}
  In older literature, the term $\infty$-category is also sometimes
  used for what we have called $(\infty,\infty)$- or
  $\omega$-categories. To avoid confusion it might therefore be
  preferable to only use the term ``$(\infty,1)$-category'', but this
  quickly gets rather cumbersome when these are the main objects being
  discussed.
\end{warning}

At first it sounds rather counterintuitive that
$(\infty,1)$-categories should be easier to define than $n$-categories
for finite $n$. The reason for this is that there is an alternative
approach to higher categories built on \emph{homotopy
  theory}. Historically, homotopy theory began as the study of
properties of topological spaces that are invariant under (weak)
homotopy equivalence. The connection to higher categories arises
through Grothendieck's \emph{Homotopy Hypothesis}, which asserts that
the homotopy-invariant information contained in a topological space is
completely captured by an associated family of higher groupoids; we
will discuss this in more detail in \S\ref{subsec:hohyp}.

Turning this idea on its head then leads to a homotopical approach to
higher categories: We can take the homotopical definition of
\igpds{} as a starting point for defining \icats{}, and then
develop other types of higher categories far more easily within the
setting of \icats{}. A key advantage of this approach is that it
allows us to avoid working combinatorially with coherence data for
associativity, as it is in a sense ``hidden away'' in the homotopy
theory of spaces. We will discuss this in more detail in
\S\ref{sec:icat}, where we consider homotopical definitions of
\icats{}, and in \S\ref{sec:incats}, where we discuss
how $(\infty,n)$-categories for $n > 1$ can be defined within the
setting of \icats{}.

\begin{remark}
  Although strict $n$-categories do not suffice for the purposes we
  are interested in here, they do have interesting connections to
  rewriting algorithms in computer science; see for instance the book
  \cite{Polygraphs} for an introduction to this topic. Strict
  $\omega$-categories also admit rather surprising algebraic
  descriptions, including that of Steiner~\cite{Steiner} in terms of
  chain complexes.
\end{remark}

\subsection{Homotopy types and fundamental $n$-groupoids}\label{subsec:htpy}
Before we can describe the Homotopy Hypothesis, we first need to
recall some ideas from algebraic topology.

\begin{defn}
  If $X$ and $Y$ are topological spaces and $f,g \colon X \to Y$ are
  continuous maps, a \emph{homotopy} from $f$ to $g$ is a continuous
  map $h \colon X \times I \to Y$ (where $I$ denotes the closed
  interval $[0,1]$) that restricts to $f$ and $g$ when the second
  coordinate is $0$ and $1$, respectively; a homotopy is thus a
  continuously varying family of maps $h_{t} = h(\blank,t)$
  interpolating between $f=h_{0}$ and $g=h_{1}$. We say that the maps
  $f$ and $g$ are \emph{homotopic} if there exists a homotopy between
  them; this is an equivalence relation on the set of continuous maps
  from $X$ to $Y$.
\end{defn}

\begin{variant}
  We can similarly consider pointed versions of these notions: A
  \emph{pointed topological space} is a pair $(X,x)$ consisting of a
  topological space $X$ and a \emph{base point} $x \in X$; a
  continuous map between pointed spaces is called \emph{pointed} if it
  preserves the given base points. For pointed continuous maps
  $f,g \colon (X,x) \to (Y,y)$ we say that a homotopy
  $h \colon X \times I \to Y$ is \emph{pointed} if $h$ takes all of
  $\{x\} \times I$ to the base point $y$ (so each of the maps $h_{t}$
  is pointed for $t \in I$). This is again an equivalence relation on
  the set of pointed continuous maps.
\end{variant}

\begin{defn}
  For a pointed topological space $(X,x)$ we define the $n$th
  \emph{homotopy group} $\pi_{n}(X,x)$ to be the set of equivalence
  classes of pointed continuous maps $(S^{n},*) \to (X,x)$ under
  pointed homotopies (where $S^{n}$ is the $n$-dimensional sphere with
  some base point $*$). Then $\pi_{1}(X,x)$ is the \emph{fundamental
    group} of loops under concatenation, while for $n > 1$ the set
  $\pi_{n}(X,x)$ has a natural abelian group structure. We also write
  $\pi_{0}(X)$ for the set of path components of $X$.
\end{defn}

\begin{defn}\label{ex:fundgpd}
  Let $X$ be a topological space. The fundamental groups of $X$ (which
  depend on a choice of base point) can be combined into the so-called
  \emph{fundamental groupoid} $\pi_{\leq 1}X$ of $X$. This category
  has the points of $X$ as its objects, and a morphism from $p$ to $q$
  is given by a homotopy class of paths in $X$ from $p$ to $q$;
  composition is given by concatenating paths. The fundamental
  groupoid is indeed a groupoid since traversing a path in the
  opposite direction gives its inverse.
\end{defn}

From the space $X$ we should also be able to extract an important
family of \emph{higher} groupoids, namely its \emph{fundamental $n$-groupoid}
$\pi_{\leq n}(X)$ for any $n \geq 0$, which also incorporates
information about all the higher
homotopy groups $\pi_{i}(X,x)$ in dimensions $i \leq n$. The
$n$-groupoid $\pi_{\leq n}(X)$ should be an $n$-category where
\begin{itemize}
\item the objects are the points of $X$, 
\item the morphisms are paths in $X$, \ie{} continuous maps $I \to X$
\item the $2$-morphisms are
homotopies between paths, \ie{} continuous maps\footnote{More
  precisely, we should impose a constancy condition on the components
  of $I \times \{0,1\}$ so that the source and target are points in $X$.} $I^{\times 2} \to X$
\item the $3$-morphisms are homotopies between
homotopies, and so on up to $n$-morphisms being equivalence classes of
$n$-dimensional homotopies in $X$ (\ie{} continuous maps
$I^{\times n} \to X$ satisfying certain constancy conditions).
\end{itemize}
Here we can also imagine that we keep going forever (and never take
homotopy classes) to obtain the \emph{fundamental \igpd{}}
$\pi_{\leq \infty}(X)$.

The homotopy hypothesis characterizes the information about $X$ that
should be contained in its fundamental $n$-groupoid $\pi_{\leq
  n}(X)$. To state this, we need some further terminology:
\begin{itemize}
\item A continuous map $f \colon X \to Y$ is a \emph{homotopy
    equivalence} if there exists a continous map $g \colon Y \to X$
  and homotopies between $gf$ and $\id_{X}$ and between $fg$ and
  $\id_{Y}$.
\item More generally, the map $f$ is a \emph{weak
  homotopy equivalence} if it induces an isomorpism
$\pi_{0}(X) \isoto \pi_{0}(Y)$ and isomorphisms
$\pi_{n}(X,x) \to \pi_{n}(Y,f(x))$ for all $n \geq 1$ and $x \in
X$.
\item We say that two spaces have the same \emph{homotopy type} if
  they are in the same equivalence class under weak homotopy
  equivalence.
\end{itemize}
Note that (by a theorem of Whitehead) a weak homotopy equivalence
between CW-complexes is actually a homotopy equivalence, so to
describe homotopy types we can either work with general topological
spaces and weak homotopy equivalences, or with nice spaces and
homotopy equivalences. All the invariants of algebraic topology are
only sensitive to the homotopy type of a topological space.

\begin{defn}
  A topological space $X$ is an \emph{$n$-type} if its homotopy groups
  $\pi_{k}(X,x)$ vanish for $k > n$ for all $x \in X$.
\end{defn}

For a topological space $X$ we can construct an $n$-type (as a
homotopy type) by ``killing'' the homotopy groups above level $n$;
this \emph{$n$-type} $\tau_{\leq n}X$ then admits a map
$X \to \tau_{\leq n}X$ that induces isomorphisms on homotopy groups in
degree $\leq n$.

\subsection{The Homotopy Hypothesis}\label{subsec:hohyp}
For small values of $n$, classical results in algebraic topology show
that there is a close relationship between $n$-groupoids and
$n$-types:
\begin{itemize}
\item A topological space is a $0$-type precisely when it is weakly
  homotopy equivalent to a set with the discrete topology. Thus we can
  identify $0$-types with sets, \ie{} $0$-groupoids.
\item Given a group $G$, the \emph{Eilenberg--MacLane space}
  \cite{EML} (or
  \emph{classifying space}) $BG$ is, up to weak homotopy
  equivalence, the unique connected space with $\pi_{1}(BG) \cong G$
  and $\pi_{n}(BG) = 0$ for $n > 1$. Moreover, this construction gives
  an equivalence of categories between groups and pointed connected
  $1$-types (with homotopy classes of maps). Since groupoids and
  general $1$-types are disjoint unions of groups and connected $1$-types,
  respectively, this correspondence extends to an equivalence between
  the $(2,1)$-categories of $1$-types and groupoids.
\item MacLane and Whitehead~\cite{MacLaneWhitehead} showed that
  connected $2$-types can be described algebraically by \emph{crossed
    modules}, which are also equivalent to strict $2$-groupoids with a
  single object \cite{BrownSpencer}. By taking disjoint unions, this
  again extends to an equivalence between $2$-types and (strict)
  $2$-groupoids.
\end{itemize}
In his unpublished manuscript \emph{Pursuing stacks}~\cite{GrothendieckStacks}, Grothendieck conjectured that this
relationship should extend to arbitrary values of $n$:
\begin{conjecture}[Grothendieck's Homotopy Hypothesis]
  There is an equivalence between $n$-groupoids and $n$-types, such
  that the fundamental $n$-groupoid $\pi_{\leq n}X$ corresponds to the
  $n$-type of the space $X$. Moreover, if we let $n$ go to infinity,
  there is an equivalence between (arbitrary) homotopy types and
  $\infty$-groupoids, such that the homotopy type of $X$ corresponds
  to the fundamental $\infty$-groupoid $\pi_{\leq \infty}X$.
\end{conjecture}

\begin{remark}
  More recently, it has also been shown that $3$-types are equivalent to
  weak $3$-groupoids; see \cite{Berger3Type} for a proof, though the
  result is originally due to the unpublished thesis of O.~Leroy. This
  equivalence also gives a concrete explanation for why strict
  $3$-groupoids are not sufficient to describe $3$-types: the strict
  $3$-groupoids only model those $3$-types that split as a product of a
  $2$-type and an Eilenberg-MacLane space $K(G,3)$
  \cite{Berger3Type}*{Corollary 3.4}; for instance, the $3$-type of the
  $2$-dimensional sphere does not correspond to a strict
  $3$-groupoid.\footnote{Simpson~\cite{Simpson3Type} has proved that
    there is also no way to realize all $3$-types by strict $3$-groupoids
    in a more general sense. In a sense the insufficiency of strict
    $3$-groupoids also goes back to work of Brown and
    Higgins~\cite{BrownHiggins}: they show that strict
    $\infty$-groupoids are equivalent to the algebraic structure of
    \emph{crossed complexes}, which are known not to model all $n$-types for $n > 2$.}
\end{remark}

Grothendieck originally formulated the Homotopy Hypothesis as a
conjecture about a hypothetical definition of $\infty$-groupoids. The
basic idea of the homotopical approach to higher categories is to
instead turn this relation on its head: we use homotopy types as a
\emph{definition} of $\infty$-groupoids, and then build more complex
notions of higher categories on top of them. This leads to a theory of
higher categories where we can avoid working with the complicated
coherence data we mentioned earlier: instead of making explicit
choices of compositions (which we then have to specify coherence data
for), we can instead merely assume that the space of possible
composites is contractible, and similarly for iterated
compositions.\footnote{We will see some ways of making this more
  precise below in \S\ref{sec:icat}.} This approach to higher
categories turns out to be much easier to work with in practice, both
for developing the general theory and for defining and working with
specific examples of higher categories.

\begin{remark}
  A classical result that partially illustrates our claim that a
  topological space ``hides'' the complex algebraic structure of an
  \igpd{} is the recognition principle for iterated loop spaces of
  May~\cite{May} and Boardman--Vogt~\cite{BoardmanVogt}: This says
  that there is an equivalence (of homotopy theories) between pointed
  $(n-1)$-connected spaces (meaning $X$ such that $\pi_{i}X = 0$ for
  $i < n$) and \emph{$E_{n}$-algebras} in topological spaces, which is
  a homotopy-coherent multiplicative structure with $k$-ary operations
  indexed by the configuration space of $k$ points in
  $\mathbb{R}^{n}$. The $E_{n}$-algebra corresponding to a space $X$
  is given by its $n$-fold loop space $\Omega^{n}X$ (or equivalently
  the space of pointed maps from $S^{n}$ to $X$). In fact, we can
  identify this $E_{n}$-algebra as precisely the structure arising
  from composition operations on the space of automorphisms of an
  identity $(n-1)$-morphism in $X$ when we view it as an \igpd{}.
\end{remark}

To give a more concrete idea of how we may use the Homotopy Hypothesis
to define higher categories, we might hope that, just as any weak
$2$-category is equivalent to a strict one, we can get away with one
level of strict associativity also here. As a first guess, we can then
simply take \emph{topological categories} as a model of \icats{}. Here
by a topological category $\uC$ we mean a category \emph{enriched} in
topological spaces, so that $\uC$ has a set of objects, and for all
objects $x,y$ a topological space $\uC(x,y)$ of morphisms, such that
the composition maps
\[ \uC(x,y) \times \uC(y,z) \to \uC(x,z)\] are continuous. If we
consider topological categories up to an appropriate notion of weak
equivalence\footnote{Cf.~\cref{defn:DKeq}.}, this does in fact turn
out to be a correct way to define \icats{}, though it has some
important drawbacks (and it certainly doesn't extend to a definition
of $(\infty,n)$-categories for $n > 1$).

In the next section we will discuss some better-behaved definitions of
\icats{}, and also attempt to explain why \icats{} are the correct
language for working with mathematical objects up to some notion of
``equivalence'' that is weaker than isomorphism.

\subsection{Further reading}
A good complement to the present article is Antol\'in Camarena's
survey \cite{Omar}, which introduces higher categories, and
especially $(\infty,1)$-categories, from a similar viewpoint to ours,
but with more focus on their applications.

Baez's paper \cite{BaezNcat} gives a more extensive introduction to
the informal idea of an $n$-category than we gave here. Cheng and
Lauda's text \cite{ChengLauda} is a very readable introduction to
several of the early approaches to $n$-categories, while Leinster's
article \cite{LeinsterSurvey} gives a brisk survey of about 10 such
definitions; his book \cite{LeinsterBook} gives a detailed
description of generalized operads and their use in defining higher categories. For a detailed discussion of equivalences in and among
higher categories, see the article \cite{ORequiv} of Ozornova and
Rovelli.

The first chapters of Simpson's book \cite{Simpson} give a good
introduction to the homotopy hypothesis and the need for weak
$n$-groupoids. The lecture notes \cite{BaezShulman} by Baez and
Shulman on $n$-categories and cohomology also give a nice introduction
to connections between topology, algebra, and higher categories in low
dimensions.

See Lack's paper \cite{Lack} for an introductory survey of
$2$-categories and Johnson and Yau's book \cite{JohnsonYau} for a
textbook treatment. The reader not intimidated by the idea of weak
$3$-categories can consult the book \cite{Gurski} by Gurski.

For more on the connection between higher categories and TQFTs, good
starting points are the original article \cite{BaezDolanTQFT} of Baez
and Dolan and Freed's expository article \cite{FreedCob}. See also
\cite{BaezLauda} for an interesting historical discussion of
connections between physics and higher category theory by Baez and
Lauda.

\section{$\infty$-categories}\label{sec:icat}
In this section we give an introduction to the homotopical theory of
\icats{}.\footnote{Recall that we always use this as an abbreviation
  for $(\infty,1)$-categories, rather than
  $(\infty,\infty)$-categories.} In recent years, this has become an
important tool in many areas of mathematics, including algebraic
topology, algebraic geometry, and representation theory, and more
generally wherever ``derived'' or ``homotopical'' structures are
employed.  What these applications have in common is that we want to
consider certain objects as equivalent in a weaker sense than just
being isomorphic, and work with constructions that are invariant under
this notion of equivalence; key examples include (weak) homotopy
equivalences between topological spaces, quasi-isomorphisms between
chain complexes, and equivalences of categories. Here we will first
explain in \S\ref{subsec:loc} why ordinary categories are not
sufficient in these situations.  We will then describe the two main
approaches to \icats{}, namely quasicategories and Segal spaces, in
\S\ref{subsec:qcat} and \S\ref{subsec:segsp}, respectively, after
reviewing simplicial sets in \S\ref{subsec:simp}, where we also
discuss simplicial categories as a first definition of \icats{}. Along
the way, we will attempt to motivate why \icats{} are a good language
for working with structures ``up to weak equivalence'', by explaining
how they arise as \emph{localizations} of ordinary categories.

\subsection{Localizations of categories}\label{subsec:loc}
Let us start by taking a quick look at the classical notion of
localizations of categories, which is a way of making weak
equivalences into isomorphisms in an ordinary category.  A
\emph{relative category} is a pair $(\uC, W)$ consisting of a category
$\uC$ and a collection $W$ of morphisms that we think of as ``weak
equivalences''. Given this data, we can always formally invert $W$:
there is a universal functor $\uC \to \uC[W^{-1}]$ that takes the weak
equivalences to isomorphisms.\footnote{This means that a functor
  $\uC \to \uD$ factors (uniquely) through $\uC[W^{-1}]$ \IFF{} it
  takes the morphisms in $W$ to isomorphisms in $\uD$.}

This construction gives, for example, the homotopy category of spaces
if we take $W$ to be the weak homotopy equivalences in the category of
topological spaces (or the homotopy equivalences if we restrict the
objects to be CW-complexes), and the derived category of a ring $R$ if
we take chain complexes of $R$-modules with quasi-isomorphisms.

An abstract definition of $\uC[W^{-1}]$ is as the pushout
$\uC \amalg_{\uW} \uW_{\txt{gpd}}$ where $\uW$ is the subcategory of
$\uC$ containing the morphisms in $W$ and $\uW_{\txt{gpd}}$ is the
groupoid obtained by formally inverting all morphisms in $\uW$ (\ie{},
$(\blank)_{\txt{gpd}}$ is the left adjoint to the inclusion of
groupoids in categories). One can also define it as a category where
the objects are those of $\uC$ and the morphisms from $x$ to $y$ are given by zig-zags
\[ x \from w \to u \from \cdots \to y,\]
where the backwards maps lie in $W$, quotiented by a certain equivalence
relation; this construction is due to Gabriel and
Zisman~\cite{GabrielZisman}.

\begin{remark}
  In general we have very little control over a localization
  $\uC[W^{-1}]$, but in many examples we can equip $\uC$ with
  additional structure that gives us more information about the
  localization. In particular, if we can equip $\uC$ with a
  \emph{model structure} \cite{QuillenHtAlg}, meaning certain classes
  of maps called ``fibrations'' and ``cofibrations'' that interact
  appropriately with the weak equivalences, then we can describe the
  set of morphisms $\Hom_{\uC[W^{-1}]}(x,y)$ as a quotient of the set
  $\Hom_{\uC}(x',y')$ where $x'$ and $y'$ are ``good'' objects weakly
  equivalent to $x$ and $y$, and we mod out by an appropriate notion of
  ``homotopies''. For example, we can describe the localization of
  topological spaces at weak homotopy equivalences as taking homotopy
  classes of maps between CW-complexes, or the derived category of a
  ring as taking chain homotopy classes of maps between complexes of
  projective modules.
\end{remark}

Unfortunately, experience shows that these ``homotopy categories''
lose a lot of important information: there are many constructions we can
make in $\uC$ and prove are invariant under weak equivalences, but
that we can't carry out just working in $\uC[W^{-1}]$.\footnote{In the
  case of derived categories and similar ``linear'' examples, we can
  however equip the homotopy category with the structure of a
  \emph{triangulated category}, which retains a lot more of this
  information.} 
\begin{ex}
  Colimits of topological spaces are not invariant under homotopy
  equivalence. For example, the $n$-dimensional sphere $S^{n}$ can be
  described as a pushout $D^{n} \amalg_{S^{n-1}} D^{n}$ of two copies
  of the $n$-dimensional unit disc $D^{n}$ (the two hemispheres of
  $S^{n}$) along its boundary $S^{n-1}$ (the ``equator''). Here
  $D^{n}$ is contractible, but the pushout
  $\mathrm{pt} \amalg_{S^{n-1}} \mathrm{pt}$ is just a point, and so
  not homotopy equivalent to $S^{n}$. Here it seems reasonable to feel
  that $S^{n}$ is in some sense the ``homotopically correct'' pushout
  (and this can be made precise as the \emph{homotopy pushout} in the
  sense of model categories), but it is \emph{not} a pushout in the
  homotopy category. In general, there exists a homotopy-invariant
  version of (co)limits of topological spaces (called \emph{homotopy
    (co)limits}), which we cannot see if we only look at the homotopy
  category.
\end{ex}

Instead of working in the homotopy category $\uC[W^{-1}]$, we can
instead try to work in $\uC$, but only make constructions that are
invariant under $W$. If we can add certain additional data to $\uC$,
such as a model category structure, we have a useful toolkit for
working with $\uC$ in this way. However, if we want to consider
objects of $\uC$ with additional structure that should be invariant
under weak equivalences, such as associative or commutative algebra
structures, model category theory becomes increasingly delicate, and
it is not of much help when we want to study interactions between
several categories that each have a notion of weak equivalence. In
these situations, it is usually much more pleasant to work in the
setting of \icats{}.

\subsection{Simplicial sets and simplicial categories}\label{subsec:simp}
Instead of using topological spaces as \igpds{}, the definitions of
\icats{} we will consider instead use \emph{simplicial sets}, which serve
as a more ``combinatorial'' description of homotopy types of
topological spaces. Here we will first recall the basic definitions of
these objects and some related notions, and then briefly consider
\emph{simplicial categories} as a first approach to defining \icats{}.

\begin{defn}
  The \emph{simplex category} $\simp$ is the category whose objects
  are the non-empty ordered sets $[n] := \{0 < 1 < \cdots < n\}$, and
  whose morphisms are the order-preserving functions. We also
  introduce notation for some basic morphisms in $\simp$:
  \begin{itemize}
  \item The $i$th \emph{coface map} $d_{i} \colon [n-1] \to [n]$ is the
    inclusion that omits $i \in [n]$, that is
    \[ d_{i}(t) =
      \begin{cases}
        t, & t < i,\\
        t+1, & t \geq i.
      \end{cases}
    \]
  \item The $i$th \emph{codegeneracy map} $s_{i} \colon [n+1] \to [n]$
    is the surjective map that repeats the value $i$, that is
    \[ s_{i}(t) =
      \begin{cases}
        t, & t \leq i,\\
        t-1, & t > i.
      \end{cases}
    \]
  \end{itemize}
\end{defn}

\begin{defn}
  A \emph{simplicial set} is a functor $\Dop \to \Set$; we write
  \[ \sSet := \Fun(\Dop, \Set)\] for the category of these. For a
  simplicial set $X$ we write $X_{n} := X([n])$ and call its elements
  the \emph{$n$-simplices} of $X$; we also write $d_{i} \colon X_{n}
  \to X_{n-1}$ and $s_{i} \colon X_{n} \to X_{n+1}$ for the face and
  degeneracy maps that are the images of the morphisms $d_{i}$ and
  $s_{i}$ in $\simp$. We denote by $\Delta^{n}$ the
  simplicial set represented by $[n] \in \simp$, so that
  $(\Delta^{n})_{k} = \Hom_{\simp}([k], [n])$; by the Yoneda lemma,
  for any simplicial set $X$ we then have a natural isomorphism $X_{n} \cong
  \Hom_{\sSet}(\Delta^{n}, X)$ so we can think of an $n$-simplex as a
  map $\Delta^{n} \to X$.
\end{defn}

\begin{defn}
  We can define a functor $\simp \to \Top$ that takes $[n]$ to the
  \emph{topological $n$-simplex}
  \[ |\Delta^{n}| := \{(x_{0},\ldots,x_{n}) \in \RR^{n+1} : 0 \leq x_{i}
    \leq 1, \, \sum x_{i} = 1\},\]
  with the subspace topology from $\RR^{n+1}$,
  and a morphism $\phi \colon [n] \to [m]$ to the continuous map given
  by
  \[ |\phi|(\mathbf{x})_{i} = \sum_{j \in \phi^{-1}(i)} x_{j}.\]
  This extends canonically to a colimit-preserving functor
  \[ |\blank| \colon \sSet \to \Top, \] the \emph{geometric
    realization}. The geometric realization has a right adjoint
  $\Sing$, the \emph{singular simplicial set} functor, with
  \[ \Sing(X)_{n} = \Hom_{\Top}(|\Delta^{n}|, X).\]
\end{defn}

We can think of a simplicial set as a ``blueprint'' for building a
topological space out of simplices through the geometric
realization. This construction also allows us to import the notion of
weak homotopy equivalences from topological spaces: A morphism of
simplicial sets $X \to Y$ is a \emph{weak equivalence} if the induced
continuous map on geometric realizations $|X| \to |Y|$ is a (weak)
homotopy equivalence. For any topological space $X$, the counit map
$|\Sing(X)| \to X$ is a weak homotopy equivalence, so that any
homotopy type can be modelled by a simplicial set. Furthermore, the
adjunction $|\blank| \dashv \Sing$ actually induces an equivalence of homotopy
categories
\[ \Top[W^{-1}] \simeq \sSet[W^{-1}],\] where $W$ denotes the weak
homotopy equivalences in $\Top$ or the weak equivalences in $\sSet$,
as appropriate.\footnote{In fact, these homotopy theories are equivalent
  in a much stronger sense: the adjunction is a \emph{Quillen
    equivalence} of model categories, which in particular means that
  they describe the same \icat{} after the \icatl{} localization we
  will discuss below.}

It is also possible to develop the homotopy theory of spaces entirely
within simplicial sets, without making reference to topological
spaces. In this theory, a special role is played by the simplicial
sets called \emph{Kan complexes}:
\begin{defn}
  For $0 \leq i \leq n$, the simplicial set $\Lambda^{n}_{i}$ (the
  $i$th \emph{horn} of $\Delta^{n}$) is the subobject of $\Delta^{n}$
  obtained by removing the interior and the $(n-1)$-dimensional face
  that does not contain the $i$th vertex. A simplicial set $X$ is a
  \emph{Kan complex} if for any map $\Lambda^{n}_{i} \to X$, we can
  find a (not necessarily unique) $n$-simplex of $X$ that extends
  it. In other words, there exists some map that fills in the diagram
  \[
    \begin{tikzcd}
      \Lambda^{n}_{i} \arrow{r} \arrow{d} & X \\
      \Delta^{n}. \arrow[dotted]{ur}
    \end{tikzcd}
  \]
\end{defn}
The Kan complexes are the most ``space-like''
objects in simplicial sets. In particular, the singular simplicial set of any topological space is a Kan complex,
and for any simplicial set $X$ we can find a weakly equivalent Kan
complex --- for example, the adjunction unit $X \to \Sing(|X|)$ is
always a weak equivalence. For a Kan complex we can, for instance,
define its homotopy groups without referring to topological spaces,
and we think of the Kan complexes as providing an alternative model
for \igpds{}.

We can now turn to simplicial categories:
\begin{defn}
  A \emph{simplicial category} is a category enriched in simplicial
  sets. In other words, a simplicial category $\uC$ consists of a set
  of objects, and for each pair of objects $x,y \in \uC$ a simplicial
  set $\uC(x,y)$ of morphisms, together with identities $\id_{x} \in
  \uC(x,x)_{0}$ and composition maps
  \[ \uC(x,y) \times \uC(y,z)\to \uC(x,z)\] that are strictly
  associative and unital; we write $\sCat$ for the category of
  simplicial categories. We define the \emph{homotopy category}
  $h\uC$ of a simplicial category $\uC$ to be the ordinary category
  obtained by taking path-components of the simplicial sets of
  maps. That is,
  \[ h\uC(x,y) = \pi_{0}(\uC(x,y)),\] where for a simplical set $X$,
  the set $\pi_{0}X$ is the quotient of $X_{0}$ by the
  relation\footnote{If $X$ is a Kan complex, then this is an
    equivalence relation.} of ``being the end points of an edge in
  $X_{1}$''
\end{defn}

\begin{remark}
  We can fully faithfully embed sets in simplicial sets as the constant
  presheaves on $\simp$. Similarly, we can view categories as those
  simplicial categories where the simplicial sets of maps are all
  constant.
\end{remark}

Given that Kan complexes can be thought of as \igpds{}, we should
expect that simplicial categories provide a model for \icats{}; at
least this should be true for the simplicial categories whose Hom's
are Kan complexes, which we call \emph{fibrant}. In order to view
simplicial categories as \icats{} we need to
consider them up to a natural notion of weak
equivalences, however:
\begin{defn}\label{defn:DKeq}
  We say a functor of simplicial categories $f \colon \uC \to \uD$ is
  a \emph{weak} (or \emph{Dwyer--Kan}) \emph{equivalence} if
  \begin{enumerate}[(1)]
  \item $f$ is \emph{weakly fully faithful}, \ie{} for $x,y \in \uC$
    the map $\uC(x,y) \to \uD(fx,fy)$ is a
    weak equivalence of simplicial sets,
  \item $f$ is \emph{essentially surjective up to homotopy}, \ie{} the
    induced functor of ordinary categories $h\uC \to h\uD$ is
    essentially surjective.
  \end{enumerate}
Any simplicial category is then weakly equivalent to one that is
fibrant.
\end{defn}
Let us consider some examples:
\begin{ex}
  We can enhance $\sSet$ to a simplicial category by taking the
  simplicial set of maps from $X$ to $Y$ to be the \emph{internal Hom}
  $Y^{X}$ in $\sSet$; this is the simplicial set $\Hom_{\sSet}(X
  \times \Delta^{\bullet}, Y)$. Let $\Kan$ denote the full
  simplicial subcategory of this where the objects are Kan complexes;
  this is a fibrant simplicial category, since the internal Hom to a
  Kan complex is again a Kan complex --- we can think of
  $\Kan$ as an incarnation of the \icat{} of \igpds{}.
\end{ex}

\begin{ex}
  Similarly, if $\uC$ is a \emph{simplicial model category}, meaning a
  model category that admits a compatible tensoring with simplicial
  sets, then we can enhance $\uC$ to a simplicial category
  $\uC_{\Delta}$: for objects $X,Y \in \uC$ the simplicial set of maps
  in $\uC_{\Delta}$ is $\Hom_{\uC}(X \otimes \Delta^{\bullet},
  Y)$. Restricting to ``nice'' (to be precise, fibrant and cofibrant)
  objects, we get a fibrant simplicial category
  $\uC_{\Delta}^{\mathrm{cf}}$. This turns out to be an incarnation of
  the \icat{} obtained from $\uC$ by inverting the weak equivalences
  (in a sense we will make more precise below).
\end{ex}

\begin{ex}
  Dwyer and Kan \cite{DwyerKan1,DwyerKan2} showed that for any
  relative category $(\uC, W)$, the localization $\uC[W^{-1}]$ is the
  homotopy category of a simplicial category $L_{W}\uC$, defined by
  one of several weakly equivalent constructions; thought of as an
  \icat{} this is again the \icatl{} localization of $\uC$ at the weak
  equivalences.\footnote{Making this precise as a universal property
    probably requires using a better notion of \icats{} than
    simplicial categories, but the paper \cite{DwyerKanDiag}
    essentially proves a version of it.}  If $\uC$ is a simplicial
  model category, then $L_{W}\uC$ (which is not fibrant) is weakly
  equivalent to the simplicial category $\uC_{\Delta}^{\mathrm{cf}}$
  from the previous example.  For a general model category $\uC$ we
  can still describe the simplicial sets of maps in $L_{W}\uC$ as
  $\Hom_{\uC}(X^{\bullet}, Y)$ using a more general kind of
  \emph{cosimplicial replacement} $X^{\bullet}$ of $X$. These
  comparisons are again due to Dwyer and Kan, in \cite{DwyerKan3}.
\end{ex}

If we are interested in developing a theory of \icats{}, there are
several disadvantages to working with simplicial categories. For
example, it is hard to give any description of the correct mapping
spaces between two simplicial categories, which means we can't define
an \icat{} of \icats{}. A related issue is that if we have a diagram
$F \colon \uI \to \uC$ (where $\uI$ is for simplicity an ordinary
category) together with equivalences\footnote{An equivalence in a
  simplicial category $\uC$ can be defined as a map $f$ such that
  composition with $f$ gives weak equivalences on all simplicial sets
  of maps.} $F(i) \isoto x_{i}$, then we can't necessarily find an
equivalent diagram $\uI \to \uC$ that takes $i$ to $x_{i}$. This was
one of the motivations for introducing \emph{homotopy-coherent
  diagrams} as a better-behaved notion of diagrams in a simplicial
category (first defined by Boardman and Vogt~\cite{BoardmanVogt} for
diagrams of topological spaces). Roughly speaking, the idea of a
homotopy-coherent diagram $F$ of shape $\uI$ in $\uC$ is that we
assign objects $F(i)$ and morphisms $F(i) \to F(j)$ in $\uC$ to
objects $i$ and morphisms $i \to j$ in $\uI$ as usual, but don't
require these assignments to respect composition. Instead, for
composable morphisms $i \xto{f} j \xto{g} k$ we assign an edge between
$F(g) \circ F(f)$ and $F(gf)$ in the simplicial set
$\uC(F(i),F(j))$. Similarly, given four composable morphisms we want a
pair of $2$-simplices
\[
  \begin{tikzcd}
    F(h)F(g)F(f) \ar[d, no head] \ar[dr, no head] \ar[r, no head] & F(h)F(gf) \ar[d,no head] \\
    F(hg)F(f) \ar[r, no head] & F(hgf)
  \end{tikzcd}
\]
relating (the images of) these edges, and so on for longer
compositions. This data amounts to a functor of simplicial categories
$\uI_{\txt{coh}} \to \uC$ from a ``coherent'' version of $\uI$; this
idea can be made precise (see \cref{rmk:Ccoh}), and one can define a
simplicial set of homotopy-coherent diagrams of shape $\uI$ in $\uC$
that indeed gives the right homotopy type of diagrams. Cordier and
Porter~\cite{CordierPorter} developed analogues of many of the main
concepts in category theory for homotopy-coherent diagrams, but the
theory is rather cumbersome to work with --- for example, it is not
easy to define the composition of two homotopy-coherent diagrams in
this language.

\subsection{Quasicategories}\label{subsec:qcat}
We are now ready to introduce the first of the two main models of
\icats{} in current use, namely \emph{quasicategories}, which are a
certain class of simplicial sets. To motivate the definition, let us
see how to describe ordinary categories as simplicial sets:

\begin{defn}
  Recall that any partially ordered set $(S, \leq)$ can be regarded as
  a category: we take the set of objects to be $S$, and say that there
  is a (unique) morphism $x \to y$ precisely when $x \leq y$. Applying
  this to the ordered sets $[n]$, we get a fully faithful embedding
  $\simp \hookrightarrow \Cat$. This extends uniquely to a
  colimit-preserving functor $\Ccat \colon \sSet \to \Cat$ with a
  right adjoint $\Nrv \colon \Cat \to \sSet$, the \emph{nerve}, given
  by
  \[ \Nrv(\uC)_{n} = \Hom_{\Cat}([n], \uC),\] so that $\Nrv(\uC)_{n}$
  is the set of composable strings of $n$ morphisms in $\uC$ (which we
  interpret to be simply objects for $n = 0$).
\end{defn}
The nerve functor is in fact fully faithful, and we can characterize
its essential image:
\begin{propn}
  A simplicial set $X$ is isomorphic to the nerve of a category \IFF{}
  every \emph{inner} horn $\Lambda^{n}_{i} \to X$ for
  $0 < i < n$ extends \emph{uniquely} to an $n$-simplex, \ie{} there
  is a unique way to fill in the diagram
  \[
    \begin{tikzcd}
      \Lambda^{n}_{i} \arrow{r} \arrow[hookrightarrow]{d} & X  \\
      \Delta^{n}. \arrow[dashed]{ur}
    \end{tikzcd}
  \]
  Furthermore, $X$ is isomorphic to the nerve of a groupoid \IFF{}
  there is a unique way to fill \emph{every} horn
  $\Lambda^{n}_{i} \to X$ for $0 \leq i \leq n$. \qed
\end{propn}
Note that the definition of Kan complexes, which we expect to be a
model for \igpds{}, is obtained by dropping the uniqueness requirement
in this characterization of nerves of groupoids. By analogy, we might
guess that we should obtain a definition of \icats{} by similarly
dropping the uniqueness in the characterization of categories:
\begin{defn}
  A \emph{quasicategory} is a simplicial set $X$ such that every inner
  horn $\Lambda^{n}_{i} \to X$ can be extended to an $n$-simplex.
\end{defn}

\begin{remark}
  This class of simplicial sets was first considered by Boardman and
  Vogt~\cite{BoardmanVogt}, who called them \emph{restricted Kan
    complexes}. They were first considered as a model of \icats{} by
  Joyal.
\end{remark}

Let us see how these objects relate to simplicial categories and
homotopy-coherent diagrams:
\begin{defn}
  Let $[n]_{\txt{coh}}$ denote the simplicial category with objects
  $0,\ldots,n$ and 
  \[ [n]_{\txt{coh}}(i,j) = \Nrv \uP_{ij},\] where $\uP_{ij}$ is the
  partially ordered set of subsets of $\{i,i+1,\dots,j\}$ that contain
  the end points; composition is induced by taking unions of
  subsets. This defines a functor $\simp \to \sCat$ that extends
  uniquely to a colimit preserving functor
  $\Ccoh \colon \sSet \to \sCat$. The functor $\Ccoh$ has a right
  adjoint $\Ncoh \colon \sCat \to \sSet$, given by
  $\Ncoh(\uC)_{n} = \Hom_{\sCat}([n]_{\txt{coh}}, \uC)$; this is
  called the \emph{coherent nerve}.
\end{defn}
If $\uC$ is a fibrant simplicial category, then $\Ncoh\uC$ is a
quasicategory, and the counit $\Ccoh\Ncoh\uC \to \uC$ is a weak
equivalence of simplicial categories. In fact, the adjunction
$\Ccoh \dashv \Ncoh$ induces an equivalence between the homotopy
theories of simplicial categories and quasicategories.

\begin{remark}\label{rmk:Ccoh}
  A homotopy-coherent diagram of shape $\uI$ in $\uC$, in the sense we
  described informally in the previous subsection, can be defined as a
  functor $\Ccoh \Nrv \uI \to \uC$ --- that is, the simplicial
  category $\Ccoh \Nrv \uI$ is precisely the ``coherent replacement''
  of $\uI$ we mentioned there. Using the adjunction
  $\Ccoh \dashv \Ncoh$, we now see that a homotopy-coherent diagram
  of shape $\uI$ in $\uC$ is the same thing as a map of simplicial
  sets $\Nrv \uI \to \Ncoh \uC$, \ie{} a diagram of shape $\uI$ in the
  quasicategory $\Ncoh \uC$.  Thus quasicategories are a setting where
  the basic notion of diagram is already a homotopy-coherent one.
\end{remark}

This ``built-in'' homotopy-coherence suggests that quasicategories
should be a very well-behaved description of \icats{}, which is in
indeed the case. For example, it is also easy to define the correct
quasicategory of functors between two quasicategories $X$ and $Y$:
this can be defined simply as the internal Hom in simplicial sets,
\ie{}
$\Fun(X,Y)_{\bullet} = \Hom_{\sSet}(X \times \Delta^{\bullet}, Y)$,
which is again a quasicategory. From this we can easily define a
version of the \icat{} of (small) \icats{}, which we could not do
using simplicial categories:
\begin{itemize}
\item Any quasicategory $\eC$ contains a maximal Kan complex $\eC^{\simeq}$ (corresponding
  to the \igpd{} of equivalences in an \icat{}).
\item We define $\QCat$ to be the (fibrant) simplicial
  category whose objects are (small) quasicategories, with the
  simplicial set of maps from $\eC$ to $\eD$ given by
  $\Fun(\eC, \eD)^{\simeq}$.\footnote{If we instead use the full
    simplicial set $\Fun(\eC, \eD)$ we obtain a simplicial category
    where the maps form quasicategories, which models the
    $(\infty,2)$-category of \icats{}.}
\item We define the quasicategory $\CatI$ of (small) \icats{} to be
  $\Ncoh(\QCat)$.
\end{itemize}
Similarly, we can define a quasicategory
$\Spc$ of spaces or \igpds{} as $\Ncoh(\Kan)$.

The \icatl{} analogues of many
constructions and theorems in category theory can be worked
out in the setting of quasicategories; this theory was first developed
by Joyal~\cite{JoyalUABNotes} and extended by Lurie~\cite{HTT} as well
as many other authors since. We do not have space to go into any
details, but we mention a few definitions to give some idea of
the flavour of the theory\footnote{See for instance \cite{HTT}*{\S 1.2} for a
longer survey of basic notions.}:
\begin{itemize}
\item If $\eC$ is a quasicategory, we can define the mapping space
  $\eC(x,y)$ for objects $x,y \in \eC_{0}$ as the pullback
  \[
    \begin{tikzcd}
      \eC(x,y) \arrow{r} \arrow{d} & \Fun(\Delta^{1}, \eC) \arrow{d}
      \\
      \{(x,y)\} \arrow{r} & \eC \times \eC.
    \end{tikzcd}
  \]
  (This is indeed a Kan complex.)
\item We say that an object $x$ in $\eC_{0}$ is terminal if the
  mapping space $\eC(x,y)$ is contractible for all $y \in \eC_{0}$. We
  can define a Kan complex of terminal objects in $\eC$ and prove this
  is always either empty or contractible --- that is, if a terminal
  object exists then it is \emph{unique} in the only way that makes
  sense if we are not allowed to distinguish between weakly
  equivalent Kan complexes.
\item Given a functor (\ie{} just a morphism of simplicial sets) $p
  \colon K \to \eC$, we can define a quasicategory of \emph{cones}
  over $p$ as the pullback
  \[
    \begin{tikzcd}
      \catname{Cone}(p) \arrow{r} \arrow{d} & \Fun(K \times \Delta^{1}, \eC) \arrow{d}
      \\
      \eC \times \{p\} \arrow{r} & \Fun(K, \eC) \times \Fun(K, \eC),
    \end{tikzcd}
  \]
  where we use the constant diagram functor $\eC \to \Fun(K, \eC)$
  (adjunct to the projection $\eC \times K \to \eC$).
\item A \emph{limit} of $p$ is a terminal object in
  $\catname{Cone}(p)$; if the limit exists, it is unique up to a
  contractible space of choices.  
\end{itemize}

\begin{remark}
  If the \icat{} $\eC$ is obtained as the localization of a model
  category $\uC$, then the \icatl{} notion of limit in $\eC$ agrees
  with the model-categorical notion of homotopy limit in
  $\uC$.\footnote{This is not hard to see given the fact that the
    \icat{} of diagrams $\Fun(\uK, \eC)$ is obtained as a localization
  of $\Fun(\uK, \uC)$, but this is itself highly non-trivial; the proof is
  one of the main topics of Cisinski's book \cite{CisinskiBook}.}
\end{remark}

\begin{remark}
  In the setting of quasicategories, we can also make precise the
  universal property that an \icatl{} localization should satisfy: a
  functor $L \colon \eC \to \eC'$ exhibits $\eC'$ as the localization
  of $\eC$ at a collection of morphisms $W$ \IFF{} for every
  quasicategory $\eD$, the functor
  $\Fun(\eC', \eD) \to \Fun(\eC, \eD)$ is fully faithful\footnote{A
    functor of quasicategories $F \colon \eX \to \eY$ is fully
    faithful if the induced maps of mapping spaces
    $\eX(x,y) \to \eY(Fx,Fy)$ are all weak equivalences.}, with image
  those functors $\eC \to \eD$ that take the morphisms in $W$ to
  equivalences. We can show that, if we replace the simplicial
  localization $L_{W}\uC$ of Dwyer--Kan by a weakly equivalent fibrant
  simplicial category $\uL$, then the resulting functor
  $\Nrv \uC \to \Ncoh \uL$ is indeed an \icatl{} localization in this
  sense.
\end{remark}

\begin{remark}
  In general, the theory of \icats{} (or quasicategories) works a lot
  like the familiar theory of categories --- once one has some
  familiarity with \icats{}, it is indeed often easy to guess what the
  \icatl{} analogue of a definition or theorem from category theory
  should be. Working out the proofs can be rather more of a challenge,
  however! A key difference between \icats{} and ordinary categories
  is that we can't just ``write down'' \icats{} and functors between
  them --- instead, we need to figure out from general principles why
  certain objects or functors should exist. For this, the description
  of functors to $\Spc$ and $\CatI$ in terms of various kinds of
  \emph{fibrations} of \icats{} becomes crucial; see \cite{BarwickShah}
  and \cite{MazelGee}
  for introductions to this important topic,
  which we do not have space to go into here.
\end{remark}

\subsection{Segal spaces}\label{subsec:segsp}
We now consider another important description of \icats{}, namely as
\emph{complete Segal spaces}, due to Rezk~\cite{RezkCSS}. The starting point is an alternative
characterization of nerves of categories:
\begin{propn}
  A simplicial set $X$ is isomorphic to the nerve of a category if
  for every $n$, the maps $X_{n} \to X_{1}, X_{0}$ induced by
  the inclusions $\{i-1,i\} \hookrightarrow [n]$ and
  $\{i\} \hookrightarrow [n]$, exhibit $X_{n}$ as an iterated fibre
  product.  In other words, the induced map
  \[ X_{n} \to X_{1} \times_{X_{0}} \cdots \times_{X_{0}} X_{1}\] is
  an isomorphism. \qed
\end{propn}
The analogue of a diagram of sets in \icat{} theory is a diagram of
\igpds{}. If we model these as simplicial sets and replace ordinary
limits by homotopy limits, we get the following definition:
\begin{defn}[Rezk, \cite{RezkCSS}]
  A \emph{Segal space} is a functor $X \colon \Dop \to \sSet$ such
  that for every $n$, the maps $X_{n} \to X_{1}, X_{0}$ induced by the
  inclusions $\{i-1,i\} \hookrightarrow [n]$ and
  $\{i\} \hookrightarrow [n]$, exhibit $X_{n}$ as an iterated
  \emph{homotopy} fibre
  product.  In other words, the induced map
  \[ X_{n} \to X_{1} \times^{h}_{X_{0}} \cdots \times^{h}_{X_{0}}
    X_{1}\] is a weak equivalence.
\end{defn}

An advantage of this definition is that it also makes sense
\emph{internally}: if we assume we have already set up a theory of
\icats{} (for instance using quasicategories), we can instead consider
Segal spaces as diagrams in the \icat{} of \igpds{}:
\begin{defn}
  A \emph{Segal space} is a functor $X \colon \Dop \to \Spc$ such
  that for every $n$, the maps $X_{n} \to X_{1}, X_{0}$ induced by the
  inclusions $\{i-1,i\} \hookrightarrow [n]$ and
  $\{i\} \hookrightarrow [n]$, exhibit $X_{n}$ as an iterated
  fibre
  product.  In other words, the induced map
  \[ X_{n} \to X_{1} \times_{X_{0}} \cdots \times_{X_{0}}
    X_{1}\] is an equivalence of \igpds{}. We write $\Seg_{\Dop}(\Spc)$ for
  the full subcategory of $\Fun(\Dop, \Spc)$ spanned by the Segal spaces.
\end{defn}
\begin{remark}
  This definition is ``model-independent'', in that it does not make
  any reference to a particular model of \icats{} (such as
  quasicategories or model categories), but is instead formulated
  internally to the \icat{} of \icats{}. If we accept the idea that
  \icats{} give a good language for working with objects up to weak
  notions of equivalence, we ought to apply this idea to \icats{}
  themselves, which is what working with \icats{} model-independently
  amounts to. In practice, this typically results in much cleaner
  statements and proofs than if we think of \icats{} as
  quasicategories, for instance, though one must take care that the
  constructions one uses are ultimately supported by the basic ones
  implemented in a model.
\end{remark}

The data of a Segal space $X$ models exactly the algebraic structure
we expect from a category:
\begin{itemize}
\item $X_{0}$ is the space of objects and $X_{1}$ is the space of
  morphisms, with $d_{1},d_{0} \colon X_{1} \to X_{0}$ indicating the
  source and target of each morphism. Given objects, that is points
  $x,y \in X_{0}$, we can thus define the space of morphisms from $x$
  to $y$ as the fibre
  \[
    \begin{tikzcd}
      X(x,y) \ar[r] \ar[d] & X_{1} \ar[d] \\
      \{(x,y)\} \ar[r] & X_{0} \times X_{0}.
    \end{tikzcd}
  \]
\item The degeneracy map $s_{0} \colon X_{0} \to X_{1}$ gives an
  identity morphism for each object.
\item The equivalence $X_{2} \simeq X_{1} \times_{X_{0}} X_{1}$
  identifies $X_{2}$ as the space of composable pairs of morphisms
  \[ x \xto{f} y \xto{g} z;\]
  their composition is given the inner face map $d_{1} \colon X_{2}
  \to X_{1}$.
\item In general $X_{n}$ can be described as the space of strings of
  $n$ composable morphisms. The image of the commutative square
  \[
    \begin{tikzcd}
      {[3]} & {[2]} \ar[l, "d_{1}"{swap}] \\
      {[2]} \ar[u, "d_{2}"]  & {[1]} \ar[l, "d_{1}"] \ar[u, "d_{1}"{swap}]
    \end{tikzcd}
  \]
  in $\Spc$ gives a homotopy between the two ways of composing 3 maps
  in two steps, so the composition is associative up to a specified
  homotopy.
\item The remaining data in the simplicial diagram then shows that the
  composition is homotopy-coherently associative and unital.  
\end{itemize}
Just as for ordinary categories, the right notion of equivalence
between Segal spaces is given by fully faithful and essentially
surjective maps, which we can define as follows:

\begin{defn}\label{def:ff}
  A morphism $F \colon X \to Y$ of Segal spaces (\ie{} a natural
  transformation $\Dop \times [1] \to \Spc$) is \emph{fully faithful}
  if the commutative square
  \[
    \begin{tikzcd}
      X_{1} \ar[r] \ar[d] & Y_{1} \ar[d] \\
      X_{0} \times X_{0} \ar[r] & Y_{0} \times Y_{0}
    \end{tikzcd}
  \]
  is a pullback.\footnote{This is equivalent to all the maps on fibres
    being equivalences, which is to say that the induced maps $X(x,y)
    \to Y(Fx,Fy)$ should be equivalences, as expected.}
\end{defn}

\begin{defn}
  Informally, an \emph{equivalence} in a Segal space $X$ is a morphism
  that has an inverse under composition. More formally, an equivalence can
  be defined as a map to $X$ from the quotient
  \[ J := \Delta^{3} \amalg_{\Delta^{1} \amalg \Delta^{1}}
    (\Delta^{0} \amalg \Delta^{0}),\]
  of $\Delta^{3}$ where
  we identify the composite edges $0 \to 2$ and $1 \to 3$ with points
  --- such a map to $X$ amounts to a diagram of the shape
  \[
    \begin{tikzcd}
      y \ar[rr, bend left=40, equals] \ar[r, "l"] & x \ar[r, "f"]
      \ar[rr, bend right=40, equals] & y \ar[r, "r"] & x,
    \end{tikzcd}
  \]
  specifying a map $f$ together with a left inverse $l$ and a right
  inverse $r$; given $f \colon \Delta^{1} \to X$, one can show that
  the space of extensions to $J$ is either empty or contractible.
\end{defn}

\begin{defn}\label{def:es}
  Let $X$ be a Segal space. Two objects in $X$ are \emph{equivalent}
  if they are connected by an equivalence; this defines an equivalence
  relation $\sim$ on $\pi_{0}X_{0}$. We say a morphism of Segal spaces
  $F \colon X \to Y$ is \emph{essentially surjective} if every object
  of $Y$ is equivalent to one in the image of $X$ --- in other words,
  if the induced map
\[ \pi_{0}X_{0}/\sim \,\to\, \pi_{0}Y_{0}/\sim \]
is surjective.
\end{defn}

To obtain the correct \icat{} of \icats{}, we should invert the fully
faithful and essentially surjective morphisms in
$\Seg_{\Dop}(\Spc)$. An advantage of the Segal space description is
that we can accomplish this localization simply by passing to a full
subcategory:
\begin{defn}\label{def:css}
  Let $X$ be a Segal space. We define $X_{1}^{\simeq}$ to be the
  subspace of $X_{1}$ consisting of equivalences.\footnote{This space
    can be shown to be equivalent to the space of maps
    $J \to X$.}  Since identities are always equivalences, the
  degeneracy map $X_{0} \to X_{1}$ factors through a map
  $X_{0} \to X_{1}^{\simeq}$. We say that $X$ is \emph{complete} if
  this map is an equivalence, and write $\CSeg_{\Dop}(\Spc)$ for the
  full subcategory of $\Seg_{\Dop}(\Spc)$ spanned by the complete
  Segal spaces.
\end{defn}
Thus $X$ is complete when its $0$th space is the ``correct'' space of
objects and equivalences between them in $X$.
\begin{thm}[Rezk~\cite{RezkCSS}]
  The inclusion $\CSeg_{\Dop}(\Spc) \hookrightarrow \Seg_{\Dop}(\Spc)$
  has a left adjoint, which exhibits $\CSeg_{\Dop}(\Spc)$ as the
  localization at the fully faithful and essentially surjective
  morphisms. \qed
\end{thm}

\begin{thm}[Joyal--Tierney~\cite{JoyalTierney}]
  There is an equivalence of \icats{}
  \[ \CSeg_{\Dop}(\Spc) \simeq \CatI \]
  between the \icat{} of complete Segal spaces and the \icat{} $\CatI$
  modelled by the fibrant simplicial category $\QCat$ of quasicategories. \qed
\end{thm}

\begin{remark}
  To\"en~\cite{ToenVers} has given axioms for the \icat{} of \icats{}
  that characterize it uniquely (up to contractible choice as always
  when working with \icats{}, and also up to the automorphism given by
  taking opposite categories).
\end{remark}

\begin{remark}
  While quasicategories and complete Segal spaces are the two most
  useful descriptions of \icats{}, there are several other models that
  are worth mentioning briefly:
  \begin{itemize}
  \item \emph{Segal categories} can be described as the Segal spaces
    $X$ such that $X_{0}$ is a set. They were first considered as a
    model of higher categories by Tamsamani~\cite{Tamsamani}, but today this
    model is perhaps mostly of historical interest; see Simpson's book
    \cite{Simpson} for an extensive treatment of the theory.
  \item As a variant of simplicial
    categories, we can consider \emph{internal} categories in
    simplicial sets. (Unlike an enriched category, which has a set of
    objects, an internal
    category has a simplicial set of
    objects.) Horel~\cite{HorelInternal} has shown that internal
    categories in $\sSet$ give a model for \icats{}.
  \item The idea of \emph{derivators}, which goes back to
    Grothendieck~\cite{GrothendieckStacks,GrothendieckDer} and
    Heller~\cite{Heller}, is that while the homotopy category of a
    relative category $(\uC, W)$ does not capture much of the
    homotopy-invariant information from $\uC$, if we remember the
    homotopy categories of $\Fun(\uI,\uC)$ for all small categories
    $\uI$ and the functors relating them we can capture a lot more ---
    in fact, we can recover the entire \icat{} from this data (see for
    instance \cite{Renaudin,Prederiv}).
  \item Barwick and Kan~\cite{BarwickKan} have shown that the
    relationship between \icats{} and localizations of ordinary
    categories is very close: one can in fact think of relative
    categories as a model of \icats{}.
  \item An analogue of quasicategories where the simplex category
    $\simp$ is replaced by a category of cubes has been developed by
    Doherty, Kapulkin, Lindsey, and Sattler~\cite{DKLSCube}.
  \end{itemize}
\end{remark}

\subsection{Further reading}
For readers previously unacquainted with simplicial sets we recommend
Friedman's introductory article \cite{Friedman}; a standard reference
for the homotopy theory of simplicial sets is the book
\cite{GoerssJardine} by Goerss and Jardine. Classic introductions to model categories include the article
\cite{DwyerSpa} by Dwyer and Spali\'nski and the book by
Hovey~\cite{Hovey}. Riehl's article \cite{RiehlHtC} is a more recent
introduction that focuses on the connection between model categories
and higher categories; see also her article \cite{RiehlCoh} for
further discussion of homotopy-coherent diagrams.

Groth's article \cite{Groth} is a short introductory survey of
quasicategories. For longer introductions we recommend the books by
Land~\cite{LandBook} and Cisinski~\cite{CisinskiBook}, and the lecture
notes by Rezk~\cite{RezkNotes} and Hinich~\cite{HinichNotes}. The book
by Riehl and Verity \cite{RiehlVerity} develops many aspects of the
theory of \icats{} in terms of an axiomatization of the simplicial
category of quasicategories (thought of as describing the
$(\infty,2)$-category of \icats{}) and its homotopy $2$-category; see
also their lecture notes \cite{RiehlVerityNotes} for an introduction
to this approach.  For more advanced topics the main source is still
Lurie's book \cite{HTT}, but this will likely eventually be superseded
by Lurie's online text in progress \cite{Kerodon}, which already
contains extensive discussions of many topics. The book \cite{Bergner}
by Bergner discusses the Quillen equivalences between quasicategories
and several other models in detail.

Once the basic framework of \icats{} is set up, one can develop
homotopy-coherent versions of algebra and algebraic geometry in this
setting. For more on this we suggest the introductory articles by
Gepner~\cite{GepnerIntro} and Rezk~\cite{RezkSAG}; Lurie's books
\cite{HA,SAG} are the standard references for the full details.

\section{$(\infty,n)$-categories}\label{sec:incats}
In this section we discuss some approaches to defining
$(\infty,n)$-categories also for $n > 1$. We first explain how to view
$(\infty,n)$-categories (and $(n,k)$-categories in general) as
\emph{enriched} \icats{} in \S\ref{subsec:enr}. Then we discuss two ways of
extending the definition of Segal spaces to higher dimensions: Rezk's
$\bT_{n}$-spaces in \S\ref{subsec:thetan} and Barwick's $n$-fold Segal
spaces in \S\ref{subsec:nfold}. 

\subsection{$(\infty,n)$-categories as enriched \icats{}}\label{subsec:enr}
Suppose $(\uV, \otimes, \bbone)$ is a monoidal category, meaning that
the category $\uV$ is equipped with
a tensor product $\blank \otimes \blank \colon \uV \times \uV \to \uV$
that is associative and unital with unit $\bbone$ up to natural
isomorphism. Then a category $\uC$ \emph{enriched} in $\uV$ has
\begin{itemize}
\item a set of objects,
\item for each pair $x,y$ of objects a morphism object $\uC(x,y)$ in
  $\uV$,
\item for each triple $x,y,z$ of objects a composition map
  \[ \uC(x,y) \otimes \uC(y,z) \to \uC(x,z), \]
\item for each object $x$ an ``identity morphism'' $\bbone \to
  \uC(x,x)$,
\end{itemize}
such that the composition is associative and unital.\footnote{To make
  this definition completely precise, we must use the natural associativity and
  unitality isomorphisms for $\uV$.}

Many naturally occurring categories can be viewed as enriched
categories --- above we have already discussed both topological and
simplicial categories, which are enriched in $\Top$ and $\sSet$ with
their cartesian products, while categories enriched in abelian groups
or chain complexes (often called \emph{dg-categories}) are ubiquitous
in algebraic settings.

We can
also view strict $n$-categories as enriched categories: if we start
with sets as $0$-categories, then a category enriched in $\Set$ is
just an ordinary category, and we can inductively define strict
$n$-categories to be categories enriched in strict $(n-1)$-categories.

This suggests that if we consider enrichment within the setting of
\icats{}, where the only notions of monoidal and enriched structures
that make sense are those that are fully homotopy-coherently
associative, we can similarly use this machinery to define \emph{weak}
higher categories. Such enriched \icats{} were first defined by Gepner
and the author in \cite{enr}. We will not review the details of
the definition\footnote{Roughly speaking, we can think of ordinary
  enriched categories as many-object versions of associative algebras
  (which can indeed be described as (pointed) enriched categories with
  a single object). Using the framework of \emph{operads} we can in
  particular describe them as algebras for many-object versions of the
  operad for associative algebras. This definition extends natually to
  the \icatl{} setting, using Lurie's machinery of
  \emph{$\infty$-operads} \cite{HA}.}
here; instead, we content ourselves with noting how a number of
higher-categorical structures can be defined in this framework:
\begin{itemize}
\item Viewing spaces as \igpds{} or $(\infty,0)$-categories, we can
  define $(\infty,n)$-categories as \icats{} enriched in
  $(\infty,n-1)$-categories. (For $n=1$, we do indeed have that
  \icats{} enriched in spaces are equivalent to \icats{}, most
  naturally in the form of complete Segal spaces \cite{enr}*{\S 4.4}.)
\item If we instead start with the category of sets, then the \icat{}
  of \icats{} enriched in $\Set$ turns out to be the $(2,1)$-category
  of ordinary categories, functors, and natural
  isomorphisms. Iterating, we get an $(n+1,1)$-category of weak
  $n$-categories as that of \icats{} enriched in weak
  $n$-categories.\footnote{This notion of weak $n$-category is
    compared to the iterated Segal categories of Tamsamani and Simpson
  in \cite{enrcomp}.}
\item We can also inductively define $(n,k)$-categories for $k \leq n$
  as \icats{} enriched in $(n-1,k-1)$-categories, where we define
  $(n,0)$-categories (or $n$-groupoids) as those $(n,1)$-categories
  where all morphisms are invertible. See \cite{enr}*{\S 6.1} for more
  details; in particular, one can show that the ``Homotopy
  Hypothesis'' is satisfied for these objects: the resulting \icat{} of
  $n$-groupoids is equivalent to that of $n$-types.
\end{itemize}

\begin{remark}
  It is often useful to extend the notion of $n$-type also to $n = -2$
  and $-1$, by declaring that the point is the only $-2$-type and the
  $-1$-types are the point and the empty set. This is compatible with
  many results on $n$-types; for instance, it remains true that a
  space is an $n$-type \IFF{} all its loop spaces are $(n-1)$-types
  also for $n = -1,0$. Viewing these $n$-types as $n$-groupoids or
  $(n,0)$-categories also for $n = -1$, we can inductively extend the
  definition of $(n,k)$-categories as \icats{} enriched in
  $(n-1,k-1)$-categories to the case $k = n+1$.\footnote{We can also
    allow $k = n+2$, but for $n > -2$ this produces nothing new: in a
    $(-1,1)$-category all mapping spaces are contractible, so it is
    either empty or equivalent to the point; these are both groupoids,
    so $(-1,1)$-categories are the same thing as
    $(-1)$-groupoids. After enriching, we get in general that
    $(n,n+2)$-categories are the same as $(n,n+1)$-categories.}
  For instance, in a $(0,1)$-category all the
  mapping spaces are either empty or contractible, so these are
  equivalent to partially ordered sets. See \cite{BaezShulman}*{\S 2
    and \S 5} for more on this topic.
\end{remark}

\begin{remark}\label{rmk:coind}
  Having defined the \icats{} $\Cat_{(\infty,n)}$ of
  $(\infty,n)$-categories, for example using enrichment, we also get a
  natural definition of $(\infty,\infty)$-categories: There are
  functors $\Cat_{(\infty,n)} \to \Cat_{(\infty,n-1)}$ that take an
  $(\infty,n)$-category to its underlying $(\infty,n-1)$-category (by
  forgetting the non-invertible $n$-morphisms), and we can define
  \[ \Cat_{(\infty,\infty)} := \limP_{n \rightarrow \infty} \Cat_{(\infty,n)}\]
  as the limit along these functors. This means that an
  $(\infty,\infty)$-category can be defined as a sequence
  $(\eC_{n})_{n\geq 0}$ where $\eC_{n}$ is an $(\infty,n)$-category
  with $\eC_{n-1}$ as its underlying $(\infty,n-1)$-category. We
  should note, however, that this notion of
  $(\infty,\infty)$-categories differs from the ``bottom-up'' approach,
  where we define an $(\infty,\infty)$-category as an algebraic
  structure with $i$-morphisms for all $i$, without defining
  $(\infty,n)$-categories for finite $n$ first. This is because there
  are problems with the concept of an $i$-morphism being invertible in
  the latter approach: we can only consider a \emph{coinductive}
  notion of equivalences, where we ``keep going forever'' in the
  definition from \S\ref{subsec:weak}. This notion may be appropriate
  in some contexts, but in general it has some very undesirable
  properties --- for instance, it turns out that all $i$-morphisms in
  the cobordism $(\infty,\infty)$-category are coinductive
  equivalences, so we cannot distinguish it from the \igpd{} where we
  invert everything. In particular, this means that the cobordism
  $(\infty,n)$-category is not the ``underlying''
  $(\infty,n)$-category of the cobordism $(\infty,\infty)$-category in
  the coinductive setting.  Goldthorpe~\cite{Goldthorpe} has shown
  that both these versions of $(\infty,\infty)$-categories have
  universal properties as fixed points for \icatl{} enrichment.
\end{remark}

\begin{remark}
  If $\uV$ is a monoidal category with a compatible monoidal
  structure, then we can define a homotopy theory of $\uV$-categories,
  which under certain assumptions is again a model category. In
  \cite{enrcomp} it is shown that the corresponding \icat{} is
  equivalent to that of \icats{} enriched in the \icatl{} localization
  of $\uV$.\footnote{If $\uV$ is symmetric monoidal, then there is an
    induced symmetric monoidal structure on $\uV$-categories, but this
    is generally \emph{not} compatible enough with the model structure
    to iterate the process --- thus this result does not contradict
    the difference between strict and weak $3$-categories, for
    instance.}  For example, \icats{} enriched in chain complexes are
  equivalent to dg-categories after \icatl{} localization. Another
  approach to weak enrichment in chain complexes is that of
  \emph{$A_{\infty}$-categories}, of particular importance in
  symplectic geometry; see for instance \cite{KellerAinfty} for an
  introduction.
\end{remark}

\subsection{$(\infty,n)$-categories as
  $\Theta_n$-spaces}\label{subsec:thetan}
In this section we introduce a definition of $(\infty,n)$-categories
due to Rezk~\cite{RezkThetan}, which generalizes Segal spaces to
higher dimensions. These objects will again be defined as presheaves of
\igpds{} that take certain diagrams to limits, and our first task is
to define the relevant indexing category:
\begin{defn}\label{def:wreath}
  For any category $\uC$, we define the \emph{wreath product}
  $\simp \wr \uC$ as follows: The objects are of the form
  $[n](c_{1},\ldots,c_{n})$ for $c_{i} \in \uC$, and a morphism
  $[n](c_{1},\ldots,c_{n}) \to [m](c'_{1},\ldots,c'_{n})$ consists of
  a morphism $\phi \colon [n] \to [m]$ in $\simp$ and a morphism
  $\psi_{ij} \colon c_{i} \to c'_{j}$ for $\phi(i-1)< j \leq \phi(j)$.
  We then inductively define the category $\bT_{n}$ as
  $\simp \wr \bT_{n-1}$, starting with $\bT_{0} = *$ (so
  $\bT_{1} = \simp$).
\end{defn}

\begin{remark}
  The object $[m](X_{1},\ldots,X_{n})$ in $\bT_{n}$ should be thought
  of as a strict $n$-category with objects $0,\ldots,n$ and with
  $X_{i}$ as the $(n-1)$-category of maps from $i-1$ to $i$ (and more
  generally with $X_{i} \times X_{i+1} \times \cdots \times X_{j}$ as
  the maps from $i-1$ to $j$). For example, the object
  $[4]([3],[0],[1],[2])$ in $\bT_{2}$ corresponds to 
  the $2$-category
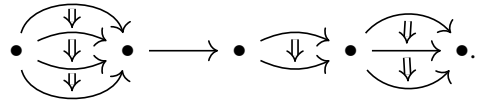
\[
\begin{tikzcd}
  \bullet 
  \ar[r, bend left=75, ""{below,name=A,inner sep=0.5pt}] 
  \ar[r, bend left=25, ""{name=B1,inner sep=0.5pt},
  ""{name=B2,below,inner sep=1pt}] 
  \ar[r, bend right=25, ""{name=C1,inner sep=1pt}, ""{below,name=C2,inner sep=0.5pt}] 
  \ar[r, bend right=75, ""{name=D,inner sep=1pt}] &
  \bullet 
  \ar[r]  &
  \bullet 
  \ar[r, bend left=25, ""{below,name=E,inner sep=1pt}] 
  \ar[r, bend right=25, ""{name=F,inner sep=1pt}] &
  \bullet 
  \ar[r, bend left=50, ""{below,name=G,inner sep=1pt}] 
  \ar[r, ""{name=H1,inner sep=1pt}, ""{name=H2,below,inner sep=1pt}] 
  \ar[r, bend right=50, ""{name=I,inner sep=1pt}] &
  \bullet.
  \arrow[from=A,to=B1,Rightarrow]
  \arrow[from=B2,to=C1,Rightarrow]
  \arrow[from=C2,to=D,Rightarrow]
  \arrow[from=E,to=F,Rightarrow]
  \arrow[from=G,to=H1,Rightarrow]
  \arrow[from=H2,to=I,Rightarrow]
\end{tikzcd}
\]  
Interpreted in this way, $\bT_{n}$ in fact defines a full subcategory
of strict $n$-categories, whose objects are \emph{free} $n$-categories
in an appropriate sense (see \cite{BergerWreath}*{Theorem 3.7}). We
think of the objects of $\bT_{n}$ as being all the basic diagram shapes
built out of $i$-morphisms for $i \leq n$ that can be composed
together in an $n$-category. For example, the free $i$-morphism
$C_{i}$ can be defined inductively in $\bT_{n}$ as
\[ C_{0} := [0](), \qquad C_{i} := [1](C_{i-1})\quad(i>0).\]
\end{remark}

\begin{remark}
  The categories $\bT_{n}$ were originally defined by
  Joyal~\cite{JoyalTheta}, who proposed a definition of higher
  categories using presheaves of sets on $\bT_{n}$. The inductive
  definition we have given is due to Berger~\cite{BergerWreath}.
\end{remark}

Now we can specify the ``Segal conditions'' we want for a presheaf on $\bT_{n}$:
\begin{defn}[Rezk, \cite{RezkThetan}]
  A \emph{Segal $\bT_{n}$-space} is a functor $X \colon \bT_{n}^{\op}
  \to \Spc$  such that
  \begin{itemize}
  \item $X([m](I_{1},\ldots,I_{m})) \isoto X([1](I_{1}))\times_{X(C_{0})}
    \cdots \times_{X(C_{0})} X([1](I_{m}))$,
  \item $X([1](\blank))$ is a Segal $\bT_{n-1}$-space.
  \end{itemize}
  We write $\Seg_{\bTnop}(\Spc)$ for the full subcategory of
  $\Fun(\bT_{n}^{\op}, \Spc)$ spanned by the Segal $\bT_{n}$-spaces.
\end{defn}

\begin{remark}
Here we think of $X(C_{i})$ as the space of $i$-morphisms for $i \leq
0 \leq n$. Unwinding the definition, it then says that the value of
$X$ at some object $I \in \bT_{n}$ is the space of composable diagrams
of $i$-morphisms in $X$ that fit together according to the shape
$I$. For instance, for $n = 2$ we have
\[ X\left(
  \begin{tikzcd}
  \bullet 
  \ar[r]  &
  \bullet 
  \ar[r, bend left=50, ""{below,name=G,inner sep=1pt}] 
  \ar[r, ""{name=H1,inner sep=1pt}, ""{name=H2,below,inner sep=1pt}] 
  \ar[r, bend right=50, ""{name=I,inner sep=1pt}] &
  \bullet
  \arrow[from=G,to=H1,Rightarrow]
  \arrow[from=H2,to=I,Rightarrow]
  \end{tikzcd}
\right) = X([2]([0],[2])) \simeq X(C_{1}) \times_{X(C_{0})}
\left(X(C_{2}) \times_{X(C_{1})} X(C_{2})\right),\] which is precisely
the space consisting of a $1$-morphism and two $2$-morphisms that fit
together in the prescribed way. The morphism $C_{2}= [1]([1]) \to
[2]([0],[2])$ given by $d_{1} \colon [1] \to [2]$ and $(s_{0} \colon
[1] \to [0], d_{1} \colon [1] \to [2])$ in $\bT_{1} = \simp$ gives a
composition map $X([2]([0],[2])) \to X(C_{2})$. In general, the
morphisms in $\bT_{n}$ precisely encode the algebraic structure of an
$(\infty,n)$-category in terms of homotopy-coherently associative and
unital compositions of $i$-morphisms.
\end{remark}

\begin{remark}
  As in the case $n = 1$, we need to invert a class of fully faithful
  and essentially surjective morphisms to obtain the correct \icat{}
  of $(\infty,n)$-categories. This can again be accomplished by
  restricting to a full subcategory of \emph{complete}
  objects\footnote{This is proved using the comparison between Segal
    $\bT_{n}$-spaces and $n$-fold Segal we discuss below.}:
\end{remark}

\begin{defn}
  We say a Segal $\bT_{n}$-space $X$ is \emph{complete} if the
  following conditions hold:
  \begin{enumerate}[(i)]
  \item The underlying Segal space of $X$, obtained by restricting
    along the fully faithful functor $\simp \hookrightarrow \bT_{n}$
    that takes $[n]$ to $[n](C_{0},\dots,C_{0})$, is complete.
  \item The Segal $\bT_{n-1}$-space $X([1](\blank))$ is complete.
  \end{enumerate}
  We write $\CSeg_{\bT_{n}}(\Spc)$ for the full subcategory of
  $\Seg_{\bT_{n}}(\Spc)$ spanned by the complete Segal $\bT_{n}$-spaces.
\end{defn}

\subsection{$(\infty,n)$-categories as $n$-fold Segal
  spaces}\label{subsec:nfold}

In this section we will consider the description of
$(\infty,n)$-categories as \emph{$n$-fold complete Segal spaces}, which
are certain presheaves of spaces on $\simp^{n}$, first introduced by
Barwick~\cite{BarwickThesis}. The starting point is the observation
that in the definition of Segal spaces we did not use any special
features of the \icat{} of spaces:
\begin{defn}
  Suppose $\eC$ is an \icat{} with finite limits. A \emph{Segal
    $\simp$-object} (or \emph{category object}) in $\eC$ is a functor
  $X \colon \Dop \to \eC$ such that the Segal maps
  \[ X_{n} \to X_{1} \times_{X_{0}} \cdots \times_{X_{0}} X_{1}\]
  are equivalences. We write $\Seg_{\Dop}(\eC)$ for the full
  subcategory of $\Fun(\Dop, \eC)$ spanned by the Segal $\simp$-objects.
\end{defn}
\begin{remark}
  This defines an \icatl{} version of \emph{internal categories}: If
  $\uC$ is a category with finite limits, then an internal category
  $X$ in $\uC$ consists of:
  \begin{itemize}
  \item objects $X_{0}, X_{1}$ (thought of as the objects and
    morphisms, respectively),
  \item maps $s,t \colon X_{1} \to X_{0}$ (assigning the source and target
    objects of the morphisms),
  \item a map $i \colon X_{0} \to X_{1}$ such that $si = \id = ti$
    (giving the identity morphisms of the objects),
  \item a map $c \colon X_{1} \times_{X_{0}} X_{1} \to X_{1}$, where the
    fibre product is over $t$ and $s$ in the two copies of $X_{1}$, such that
    \[ sc = s\circ \txt{pr}_{1}, \quad tc = t \circ \txt{pr}_{2},
      \quad
      c \circ (\id \times_{X_{0}} i) = \id = c \circ (i \times_{X_{0}}
      \id)\]
      (giving
    the composites of composable pairs of morphisms),
  \item and such that we have a commutative square
    \[
      \begin{tikzcd}
        X_{1} \times_{X_{0}} X_{1} \times_{X_{0}} X_{1} \ar[r, "c
        \times_{X_{0}} \id"] \ar[d, "\id \times_{X_{0}} c"] &  X_{1}
        \times_{X_{0}} X_{1} \ar[d, "c"] \\
         X_{1} \times_{X_{0}} X_{1} \ar[r, "c"] & X_{1},
       \end{tikzcd}
     \]
     representing the associativity of the composition.
   \end{itemize}
   Note that an internal category in $\Set$ is just an ordinary
   category.\footnote{An ordinary category is of course also a category
   \emph{enriched} in $\Set$, but in all other cases internal and
   enriched categories are very different!}
\end{remark}

\begin{ex}
  An internal category $X$ in $\Cat$ is a \emph{double category}: we
  think of the objects of $X_{0}$ as objects, the morphisms in $X_{0}$
  as \emph{vertical} morphisms, and the objects of $X_{1}$ as
  \emph{horizontal} morphisms; the morphisms in $X_{1}$ we can then
  think of as \emph{squares}:
  \[
    \begin{tikzcd}
      \bullet \ar[r] \ar[d] \ar[dr, phantom, "\Downarrow"] & \bullet \ar[d] \\
      \bullet \ar[r] & \bullet,
    \end{tikzcd}
  \]
  with the top and bottom horizontal arrows representing the source
  and target objects in $X_{1}$ and the left and right vertical arrows
  the source and target morphisms in $X_{0}$. We can compose these
  squares vertically (using the composition in $X_{1}$) and
  horizontally (using the composition functor
  $c \colon X_{1} \times_{X_{0}} X_{1} \to X_{1}$), and these
  compositions are compatible (since $c$ is a functor).
\end{ex}

We can iterate the internal category construction, and inductively
define an \emph{$n$-uple category}\footnote{These are more commonly
  called \emph{$n$-fold categories}, but to avoid confusion we reserve
  the word ``$n$-fold'' for $n$-fold Segal spaces.} as an internal
category in $(n-1)$-uple categories (starting with sets as $0$-uple
categories). 
Then an $n$-uple category has objects, $n$ different
types of morphisms, $\binom{n}{2}$ different types of square-shaped
$2$-morphisms, and in general $\binom{n}{k}$ different types of
$k$-morphisms in the form of $k$-dimensional cubes. We can also
consider the analogous structures in the \icatl{} setting:
\begin{defn}
  Suppose $\eC$ is an \icat{} with finite limits. Then a \emph{Segal
    $\simp^{n}$-object} (or \emph{$n$-uple category object}) in $\eC$
  is inductively defined to be a Segal $\simp$-object in the \icat{}
  of Segal $\simp^{n-1}$-objects in $\eC$. We write
  $\Seg_{\Dnop}(\eC)$ for the \icat{}
  $\Seg_{\simp^{\op}}(\Seg_{\simp^{n-1,\op}}(\eC))$ of Segal $\simp^{n}$-objects in
  $\eC$; this is a full subcategory of $\Fun(\Dnop, \eC)$. We refer to
  Segal $\simp^{n}$-objects in the \icat{} of spaces as \emph{Segal
    $\simp^{n}$-spaces} (or \emph{$n$-uple
  Segal spaces}).
\end{defn}

\begin{remark}
  Viewing $\CatI$ as the \icat{} of complete Segal spaces, we can
  think of the \icat{} $\Seg_{\simp^{n-1}}(\Cat)$ as a full subcategory of
  $\Seg_{\simp^{n}}(\Spc)$. In practice, $n$-uple Segal spaces often
  arise as such Segal $\simp^{n-1}$-\icats{} (or \emph{$n$-uple
    \icats{}}). Note, however, that it typically does not make sense
  to impose further completeness conditions, as the equivalences for
  the different types of morphisms are usually different.
\end{remark}

To connect these structures to $(\infty,n)$-categories, let us go back
to the setting of ordinary double categories, and observe that we can
identify strict $2$-categories as the double categories where the only
vertical morphisms are identities (or in other words, the double
categories $X$ where $X_{0}$ is a \emph{set}); in effect, we think of
the $2$-morphisms as squares of the form
\[
  \begin{tikzcd}
    \bullet \ar[r] \ar[d, equals] \ar[dr, phantom, "\Downarrow"] &
    \bullet \ar[d, equals] \\
    \bullet \ar[r] & \bullet.
  \end{tikzcd}
\]
In the same way, we can identify strict $n$-categories as the $n$-uple
categories where all but one type of $i$-morphism is trivial for each
$i$. The idea of $n$-fold Segal spaces is to similarly identify
$(\infty,n)$-categories as $n$-uple Segal spaces satisfying certain
constancy conditions:
\begin{defn}[Barwick, \cite{BarwickThesis}]
  We inductively define the \emph{$n$-fold Segal spaces} to be the
  Segal $\simp^{n}$-spaces $X \colon \simp^{n,\op} \to \Spc$ such that
  \begin{enumerate}[(i)]
  \item $X_{0} \colon \simp^{n-1,\op} \to \Spc$ is constant,
  \item $X_{1} \colon \simp^{n-1,\op} \to \Spc$ is an $(n-1)$-fold
    Segal space.
  \end{enumerate}
  If $n = 1$, any Segal space is a $1$-fold Segal space.  We write
  $\Seg_{(n)}(\Spc)$ for the full subcategory of $\Seg_{\Dnop}(\Spc)$
  spanned by the $n$-fold Segal spaces.
\end{defn}

Just as in the case $n = 1$, we need to invert a class of ``fully
faithful and essentially surjective'' morphisms to get the correct
\icat{} of $(\infty,n)$-categories, and we can do this by restricting
to a full subcategory of complete objects\footnote{This localization
  follows from the inductive presentation of complete $n$-fold Segal
  spaces in \cite{LurieGoodwillie}.}:
\begin{defn}
  A morphism $X \to Y$ of $n$-fold Segal spaces is \emph{fully
    faithful and essentially surjective} if
  \begin{enumerate}[(i)]
  \item $X_{\bullet,0,\dots,0} \to Y_{\bullet,0,\dots,0}$ is a fully
    faithful and essentially surjective morphism of Segal spaces in
    the sense of Definitions~\ref{def:ff} and \ref{def:es},
  \item $X_{1} \to Y_{1}$ is a fully faithful and essentially
    surjective morphism of $(n-1)$-fold Segal spaces.
  \end{enumerate}
\end{defn}

\begin{defn}
  We say that an $n$-fold Segal space $X$ is \emph{complete} if
  \begin{enumerate}[(i)]
  \item $X_{\bullet,0,\dots,0} \colon \simp^{\op} \to \Spc$ is a
    complete Segal space in the sense of \cref{def:css},
  \item $X_{1} \colon \simp^{n-1,\op} \to \Spc$ is a complete
    $(n-1)$-fold Segal space.
  \end{enumerate}
  We write $\CSeg_{(n)}(\Spc)$ for the full subcategory of
  $\Seg_{(n)}(\Spc)$ spanned by the complete $n$-fold Segal spaces.
\end{defn}

\begin{remark}
  The \icat{} $\CSeg_{(n)}(\Spc)$ is equivalent to the \icat{}
  $\Cat_{(\infty,n)}$ defined by iterated enrichment; this is shown in
  \cite{enrcomp}*{\S 7}.
\end{remark}

\begin{remark}
  For any category $\uC$ we can define a functor from $\simp \times
  \uC$ to the wreath product $\simp \wr \uC$ of \cref{def:wreath} by
  taking $([n], I)$ to $[n](I,\dots,I)$; iterating this definition
  gives a functor $\tau_{n} \colon \simp^{n} \to \bT_{n}$. Composition
  with this functor restricts to equivalences
  \[ \tau_{n}^{*} \colon \Seg_{\bT_{n}}(\Spc) \isoto \Seg_{(n)}(\Spc),
    \qquad \tau_{n}^{*} \colon \CSeg_{\bT_{n}}(\Spc) \isoto
    \CSeg_{(n)}(\Spc).\]
  This was first proved by Barwick and
  Schommer-Pries~\cite{BSP}; other proofs have also been given by
  Bergner and Rezk~\cite{BergnerRezk} and the author \cite{thetan}.
\end{remark}

\begin{remark}\label{rmk:exs}
  We revisit the examples of $n$-categories we sketched in
  \S\ref{subsec:exs} and give some references for their
  $(\infty,n)$-categorical versions:
  \begin{itemize}
  \item \cref{ex:catofncat} generalizes to the $(\infty,n+1)$-category
    of (small) $(\infty,n)$-categories; this is perhaps most naturally
    viewed as the self-enrichment of $\Cat_{(\infty,n)}$ that follows
    from this being a cartesian closed \icat{}. Such enrichments are
    obtained by a somewhat ad hoc construction in \cite{enr}*{\S 7};
    Hinich gives a more natural construction in \cite{HinichYoneda}.
  \item The Morita $2$-categories of algebras and bimodules from
    \cref{ex:morita} can be extended to higher dimensions as Morita
    $(\infty,n+1)$-categories of \emph{$E_{n}$-algebras} and iterated
    bimodules. See \cite{nmorita} for an algebraic construction and
    Scheimbauer's thesis~\cite{ScheimbauerThesis} for a more geometric
    version defined using factorization algebras.
  \item The higher categories of spans from \cref{ex:span}, and more
    generally $(\infty,n)$-categories of iterated spans
    in an \icat{} with finite limits, are defined in \cite{spans}.
  \item The cobordism $n$-category of \cref{ex:cob} naturally extends
    to an $(\infty,n)$-category that includes diffeomorphisms, smooth
    homotopies, and so on as invertible $i$-morphisms for $i >
    n$. Such cobordism $(\infty,n)$-categories are constructed in
    \cite{CSBord}.
  \end{itemize}
\end{remark}

\begin{remark}\label{rmk:otherinftyn}
  Let us briefly mention a few other models of $(\infty,n)$-categories
  that we do not have space to discuss in detail:
  \begin{itemize}
  \item A \emph{stratified simplicial set} is a simplicial set
    equipped with a collection of ``marked'' $n$-simplices for each $n >
    0$, which must include the degenerate ones.
    erity~\cite{VerityCompl} has proposed a definition
    of $(\infty,\infty)$-categories as stratified simplicial sets that
    satisfy a collection of lifting properties, called
    \emph{complicial sets}.\footnote{Or rather \emph{weak} complicial
      sets in Verity's paper, with complicial sets originally
      referring to the stratified simplicial sets that describe strict
      $(\infty,\infty)$-categories \cite{VerityStrict}.}
    The basic idea, which goes back to work of
    Street~\cite{StreetOriented} is that we think of $0$- and
    $1$-simplices of a simplicial set $X$ as objects and morphisms,
    just as for quasicategories, but the $2$-simplices should now be
    $2$-categorical diagrams of the form
    
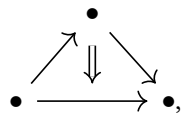
\[
      \begin{tikzcd}[column sep=small]
        {} & \bullet \ar[dr] \\
        \bullet \ar[ur] \ar[rr, ""{above,name=A}] & &
        \bullet,
        \ar[from=1-2, to=A, Rightarrow, shorten <=5pt, shorten >= 1pt]
      \end{tikzcd}
    \]
    the $3$-simplices should be tetrahedra whose faces have such
    $2$-morphisms and which contain a $3$-morphism, and so forth. The marked
    $n$-simplices then pick out those diagrams of this type where the
    $n$-morphism is invertible.\footnote{The fact that the invertible
      $n$-morphisms must be specified in this way is closely related
      to the difficulty of defining invertible morphisms in an
      $(\infty,\infty)$-category that we discussed in
      \cref{rmk:coind}; indeed, this data is in a sense superfluous if
      we only consider the complicial sets that describe
      $(\infty,n)$-categories for some finite $n$.}  If we assume
    further that all $k$-simplices are marked for $k > n$, then these
    \emph{$n$-trivial} complicial sets should be a model for
    $(\infty,n)$-categories. For $n = 1$, the $1$-trivial complicial
    sets are just quasicategories marked by their equivalences. The
    comparison is also known for $n = 2$ by results of
    Lurie~\cite{LurieGoodwillie} and Gagna, Harpaz, and
    Lanari~\cite{GHL}, but is still open for $n > 2$.
  \item Campion, Kapulkin and Maehara~\cite{CKMCube} have defined an
    analogue of complicial sets using cubical sets, which is compared to
    Verity's simplicial version in \cite{DKMCube}.
  \item One can consider analogues of quasicategories in presheaves of
    \emph{sets} on $\bT_{n}$; this definition has been worked out by
    Ara~\cite{AraTheta}, who shows that it gives a model equivalent to
    Rezk's Segal $\bT_{n}$-spaces.
  \item Barwick and Schommer-Pries~\cite{BSP} show that
    $(\infty,n)$-categories can be viewed as presheaves of spaces on
    \emph{gaunt $n$-categories}, which are the strict $n$-categories
    without any non-trivial invertible $i$-morphisms for all $i$.
  \item Barwick and Kan~\cite{BarwickKanN} show that
    $(\infty,n)$-categories can be modelled by \emph{$n$-relative
      categories}, which are defined to be certain categories equipped
    with $n+1$ wide subcategories.
  \end{itemize}
\end{remark}

\begin{remark}
  An important topic of current research is lax transformations, lax
  functors, and Gray tensor products. For $(\infty,2)$-categories
  these notions can be defined both via $2$-fold Segal spaces (viewed as
  fibrations over $\Dop$), as in \cite{GR}, and via scaled simplicial
  sets (\ie{} $2$-trivial complicial sets), as in \cite{GHLGray}; these
  two versions have recently been compared by
  Abell\'an~\cite{Abellan}. Many open questions remain, however, in
  particular regarding the equivalence of different versions of Gray
  tensor products and lax transformations for $(\infty,n)$-categories,
  (including for $n = 2$, but especially for $n > 2$); see for
  instance \cite{JFS,ORV,CampionGray}.
\end{remark}

\subsection{Further reading}
Bergner's survey \cite{BergnerSurvey} discusses several models of
$(\infty,n)$-categories from a
model-categorical perspective, while Riehl's lecture notes
\cite{RiehlCompl} give an introduction to Verity's complicial
sets. The appendix of the book of Gaitsgory and Rozenblyum \cite{GR}
discusses many constructions and results for $(\infty,2)$-categories,
and has motivated a lot of recent work on the subject.

The theory of $(\infty,n)$-categories is still under active
development, and so far not too many expository texts have appeared;
unlike in the previous sections we will therefore also point the
reader to a few interesting research papers:
\begin{itemize}
\item Barwick and Schommer-Pries~\cite{BSP} give axioms that uniquely
  characterize the \icat{} of $(\infty,n)$-categories (up to the
  automorphisms given by reversing $i$-morphisms).
\item Nuiten~\cite{Nuiten} proves a straightening theorem for
  fibrations of $(\infty,n)$-categories.
\item Important work on enriched \icats{} includes the papers of
  Hinich~\cite{HinichYoneda} on the Yoneda lemma and of
  Heine~\cite{Heine} on enrichment from (weak) module structures.
\end{itemize}

\section*{Acknowledgments}
I thank Joachim Kock for helpful comments on a draft of this article.

\end{document}